\documentclass[a4paper,reqno,11pt]{amsart}
\usepackage{amsfonts,amsmath,amsthm,amssymb,stmaryrd,mathrsfs}
\usepackage[numbers,sort&compress]{natbib}
\newtheorem{theorem}{Theorem}[section]
\newtheorem{lemma}[theorem]{Lemma}

\newtheorem{Remark}[theorem]{Remark}

\newtheorem{Corollary}[theorem]{Corollary}
\numberwithin{equation}{section}
\allowdisplaybreaks

\usepackage[left=1 in, right=1 in,top=1 in, bottom=1 in]{geometry}

\newcommand{\na}{\nabla}

\newcommand{\al}{\alpha}

\newcommand{\ga}{\gamma}

\newcommand{\la}{\lambda}

\providecommand{\norm}[1]{\left\Vert#1\right\Vert}

\providecommand{\norms}[1]{\left\vert#1\right\vert}
\def\r3{\mathbb{R}^3}

\begin{document}
\title[Compressible Euler-Maxwell system]{Decay estimates of solutions to the compressible Euler-Maxwell system in $\r3$}

\author{Zhong Tan}
\address{School of Mathematical Sciences\\
Xiamen University\\
Xiamen, Fujian 361005, China}
\email[Z. Tan]{ztan85@163.com}

\author{Yanjin Wang}
\address{School of Mathematical Sciences\\
Xiamen University\\
Xiamen, Fujian 361005, China}
\email[Y. J. Wang]{yanjin$\_$wang@xmu.edu.cn}

\author{Yong Wang}
\address{School of Mathematical Sciences\\
Xiamen University\\
Xiamen, Fujian 361005, China}
\email[Y. Wang]{wangyongxmu@163.com}



\thanks{Corresponding author: Yong Wang, wangyongxmu@163.com}

\thanks{Supported by the National Natural Science Foundation of China (Grant No. 11271305) and the Natural Science Foundation of Fujian Province of China (No. 2012J05011).}

\begin{abstract}
We study the large time behavior of solutions near a constant
equilibrium to the compressible Euler-Maxwell system in $\r3$. We
first refine a global existence theorem by assuming that the $H^3$
norm of the initial data is small, but the higher order derivatives
can be arbitrarily large. If the initial data belongs to
$\Dot{H}^{-s}$ ($0\le s<3/2$) or $\dot{B}_{2,\infty}^{-s}$ ($0<
s\le3/2$), by a regularity interpolation trick, we obtain the
various decay rates of the solution and its higher order
derivatives. As an immediate byproduct, the usual $L^p$--$L^2$
$(1\le p\le 2)$ type of the decay rates follow without requiring
that the $L^p$ norm of initial data is small.
\end{abstract}

\maketitle

\section{Introduction}
%%%%%%%%%%%%%%%%%%%%%%%%%%%%%%%%%%%%%%%%%%%%%%

The dynamics of the electrons interacting with their
self-consistent electromagnetic field can be described by the compressible Euler-Maxwell system \cite{MRS}:
\begin{equation}  \label{yiyi}
\left\{
\begin{array}{lll}
\partial_t \tilde n+{\rm div} (\tilde n \tilde u)=0,\\
\partial_t (\tilde n \tilde u)+{\rm div}(\tilde n \tilde u\otimes \tilde u) + \na p(\tilde n)=-\tilde n(\tilde{E}+\varepsilon \tilde u\times \tilde{B})- \frac{1}{\tau} \tilde n \tilde u,\\
\varepsilon\lambda^2\partial_t \tilde{E}-\na \times \tilde{B}= \varepsilon\tilde n \tilde u,\\
\varepsilon\partial_t \tilde{B}+\na \times \tilde{E}=0,\\
\lambda^2{\rm div}\tilde{E}= n_\infty  -\tilde n,\ \ {\rm div}\tilde{B}=0,\\
(\tilde n,\tilde u,\tilde{E},\tilde{B})|_{t=0}=(\tilde n_0, \tilde u_0, \tilde{E}_0, \tilde{B}_0).
\end{array}
\right.
\end{equation}
The unknown functions $\tilde n, \tilde u, \tilde{E}, \tilde{B} $ represent the electron density,
electron velocity, electric field and magnetic field, respectively. We assume the pressure $p(\tilde n)=A\tilde n^{\gamma}$
with constants $A>0$ and $\gamma\ge 1$  the adiabatic exponent. $\tau>0$ is the relaxation time. $\lambda>0$ is the Debye length, and $\varepsilon=1/c$ with $c$ the speed of light. In the motion of the fluid, due to the greater inertia the ions merely provide a constant charged
background $ n_\infty  >0$.

Despite its physical importance, due to the complexity there are only few mathematical studies on the Euler-Maxwell system. In one space dimension, Chen, Jerome and Wang \cite{CJW} proved the global existence of entropy weak solutions to the initial-boundary value problem for arbitrarily large initial data in $L^{\infty}$. Since the Euler-Maxwell system is a symmetrizable hyperbolic system, the Cauchy problem in $\r3$ has a local unique smooth solution when the initial data is smooth, see Kato \cite{K} and Jerome \cite{J} for instance. Recently, there are some results on the global existence
and the large time behavior of smooth solutions with small perturbations, see Duan \cite{D}, Ueda and Kawashima \cite{UK},  Ueda, Wang and Kawashima \cite{UWK}. For the asymptotic limits that derive simplified models starting from the Euler-Maxwell system, we refer to \cite{HP,PWG,X} for the relaxation limit; \cite{X} for the non-relativistic limit; \cite{PW1,PW2} for the quasi-neutral limit; \cite{T1,T2} for WKB asymptotics; and references therein.

The main purpose of this paper is to derive some various time decay rates of the solution as well as its spatial derivatives of any order. Meantime, we also establish a refined global existence of smooth solutions near the constant equilibrium $(n_\infty,0,0,B_\infty)$ to the compressible Euler-Maxwell system, compared with \cite{D,UWK}. We should emphasize that our results highly rely on that we consider the relaxation case. The non-relaxation case is much more difficult, we refer to \cite{GM,G,IP} for such direction. It turns out that it is more convenient to reformulate the compressible Euler-Maxwell system \eqref{yiyi} as follows.
Without loss of generality, we take the constants $\tau,\varepsilon,\lambda ,A, n_\infty$ in \eqref{yiyi} to be one.
We define
\begin{equation}
\left\{
\begin{array}{lll}
n(x,t)=\frac{2}{\gamma-1}\left\{\left[\tilde n\Big(x,\frac{t}{\sqrt{\gamma}}\Big)\right]^{\frac{\gamma-1}{2}}-1\right\},\quad u(x,t)=\frac{1}{\sqrt{\gamma}}\tilde{u}\Big(x,\frac{t}{\sqrt{\gamma}}\Big), \\
E(x,t)=\frac{1}{\sqrt{\gamma}}\tilde{E}\Big(x,\frac{t}{\sqrt{\gamma}}\Big), \qquad B(x,t)=\frac{1}{\sqrt{\gamma}}\tilde{B}\Big(x,\frac{t}{\sqrt{\gamma}}\Big)- B_\infty .
\end{array}
\right.
\end{equation}
Then the Euler-Maxwell system \eqref{yiyi} is reformulated equivalently as
\begin{equation}  \label{EM per}
\left\{
\begin{array}{lll}
\displaystyle\partial_tn+{\rm div} u=-u\cdot\na n
-\mu n{\rm div} u,   \\
\displaystyle\partial_tu+\nu  u+u\times  B_\infty +\na n+\nu  E=-u\cdot\na u
-\mu n\na n-u\times B, \\
\partial_t E-\nu  \na \times B-\nu  u= \nu  f(n) u, & & \\
\partial_t B+\nu  \na \times E=0,  \\
{\rm div} E=-\nu  f(n),\ \ {\rm div} B=0, \\
(n,u,E,B)|_{t=0}=(n_0, u_0, E_0, B_0).
\end{array}
\right.
\end{equation}
Here $\mu:=\frac{\gamma-1}{2}$, $\nu:=\frac{1}{\sqrt{\gamma}}$ and the nonlinear function $f(n)$ is defined by
\begin{equation}  \label{f}
f(n):=\left(1+ \frac{\gamma-1}{2} n\right)^{\frac{2}{\gamma-1}}-1.
\end{equation}
Notice that we have assumed $\ga>1$. If $\ga=1$, we instead define
\begin{equation}
n:=\sqrt{A}\left({\rm ln}\ \tilde n-{\rm ln}\  n_\infty  \right)=\sqrt{A} {\rm ln}\ \tilde n.
\end{equation}
In this paper, we only consider the case $\ga>1$, and the case $\ga=1$ can be treated in the same way by using the reformulation in terms of the new variables correspondingly.
\smallskip\smallskip

\noindent \textbf{Notation.} In this paper, we use $H^{s}(\mathbb{R}^{3})$, $s\in \mathbb{R}$ to denote the usual
Sobolev spaces with norm $\norm{\cdot}_{H^{s}}$ and $L^{p}(\mathbb{R}^{3})$, $1\leq p\leq
\infty $ to denote the usual $L^{p}$ spaces with norm $
\norm{\cdot}_{L^{p}}$.
$\na ^{\ell }$ with an
integer $\ell \geq 0$ stands for the usual any spatial derivatives of order $
\ell $. When $\ell <0$ or $\ell $ is not a positive integer, $\na ^{\ell }
$ stands for $\Lambda ^{\ell }$ defined by $
\Lambda^\ell f := \mathscr{F}^{-1} (|\xi|^\ell  \mathscr{F}{f}) $, where $\mathscr{F}$ is the usual Fourier transform operator and $\mathscr{F}^{-1}$ is its inverse. We use $\dot{H}
^{s}(\mathbb{R}^{3})$, $s\in \mathbb{R}$ to denote the homogeneous Sobolev
spaces on $\mathbb{R}^{3}$ with norm $\norm{\cdot}_{\dot{H}^{s}}$ defined by
$\norm{f}_{\dot{H}^s}:=\norm{\Lambda^s f}_{L^2}$. We then recall the homogeneous Besov spaces. Let
$\phi \in C^\infty_c(R^3_\xi)$ be such that $\phi(\xi) = 1$ when $|\xi| \le 1$ and $\phi(\xi) = 0$ when $|\xi| \ge 2$.  Let
$\varphi(\xi) = \phi(\xi) - \phi(2\xi)$ and $\varphi_j(\xi) = \varphi( 2^{-j}\xi)$
for ${j \in \mathbb {Z}}$. Then by the
construction,
$\sum_{j\in\mathbb {Z}}\varphi_j(\xi)=1$ if $\xi\neq 0.$
We define $\dot{\Delta}_j
f:= \mathscr{F}^{-1}(\varphi_j)* f$, then for $s\in \mathbb{R}$ we define the homogeneous Besov spaces $\dot{B}_{2,\infty}^{s}(\r3)$ with norm $\norm{\cdot}_{\dot{B}_{2,\infty}^{s}}$ defined by
\begin{equation}
\|f\|_{\dot{B}_{2,\infty}^{s}}:=\sup\limits_{j\in\mathbb{Z}}2^{sj}\norm{\dot{\Delta}_j f}_{L^{2}}.
\end{equation}

Throughout this paper  we let $C$  denote
some positive (generally large) universal constants and $\lambda$ denote  some positive (generally small) universal constants. They {\it do not} depend on either $k$ or $N$; otherwise, we will denote them by $C_k$, $C_{N}$, etc.
We will use $a \lesssim b$ if $a \le C b$, and $a\thicksim b$ means that $a\lesssim b$ and $b\lesssim a$.
We use $C_0$ to denote the constants depending on the initial data and $k,N,s$.
For simplicity, we write $\norm{(A,B)}_{X}:=\norm{A}_{X} +\norm{B}_{X}$ and $\int f:=\int_{\mathbb{R}^3}f\,dx.$

\smallskip\smallskip

For $N\ge 3$, we define the energy functional by
\begin{equation}
\mathcal{E}_N(t):=\sum_{l=0}^{N}\norm{ \na^{l}(n, u, E, B )}_{L^2}^2
\end{equation}
and the corresponding dissipation rate by
\begin{equation}\label{DN}
\mathcal{D}_N(t):=\sum_{l=0}^{N}\norm{\na^{l}(n,u)}_{L^2}^2+\sum_{l=0}^{N-1}\norm{\na^{l}E}_{L^2}^2+ \sum_{l=1}^{N-1}\norm{\na^{l} B}_{L^2}^2.
\end{equation}

As a byproduct of our analysis for deriving the decay rate of the
solution to the system \eqref{EM per}, we may first refine a global
existence theorem as stated in the following.
\begin{theorem}\label{existence}
Assume the initial data satisfy the compatible conditions
\begin{equation}\label{compatible condition}
{\rm div}  {E}_0=-\nu  f(n_0),\ \ {\rm div}  {B}_0=0.
\end{equation}
There exists a sufficiently small $\delta_0>0$ such that if $\mathcal{E}_3(0)\le \delta_0$, then there exists a unique global solution
$(n,u,E,B)(t)$ to the Euler-Maxwell system \eqref{EM per} satisfying
\begin{equation}\label{energy inequality}
\sup_{0\leq t\leq \infty }\mathcal{ E}_3(t)+\int_{0}^{\infty }
\mathcal{D}_3(\tau)\,d\tau\leq C\mathcal{ E}_3(0).
\end{equation}

Furthermore, if $\mathcal{E}_N(0)<+\infty$ for any $N\ge 3$, there exists an increasing continuous function $P_N(\cdot)$ with $P_N(0)=0$ such that the unique solution satisfies
\begin{equation}\label{energy inequality N}
\sup_{0\leq t\leq \infty }\mathcal{ E}_N(t)+\int_{0}^{\infty }
\mathcal{D}_N(\tau)\,d\tau\leq P_N\left(\mathcal{ E}_N(0)\right).
\end{equation}
\end{theorem}

The proof of Theorem \ref{existence} is inspired by the recent work
of Guo \cite{G12}. The new major difficulty here is the
regularity-loss of the electromagnetic field. We will do the refined
energy estimates stated in Lemma \ref{energy lemma}--\ref{other di},
which allow us to deduce
\begin{equation}
\frac{d}{dt} \mathcal{E} _3+\mathcal{D}_3\lesssim\sqrt{\mathcal{E} _3}\mathcal{D}_3
\end{equation}
and for $N\ge 4$,
\begin{equation}
\frac{d}{dt} {\mathcal{E}}_N+ \mathcal{D}_N
 \le  C_{N } {\mathcal{D}_{N-1}} {\mathcal{E}_N}.
\end{equation}
Then Theorem \ref{existence} follows in the fashion of \cite{G12}.

The main purpose of this paper is to derive some various decay rates
of the solution to the system \eqref{EM per} by making the stronger
assumption on the initial data.
\begin{theorem}\label{decay}
Assume that $(n,u,E,B)(t)$ is the solution to the
Euler-Maxwell system \eqref{EM per} constructed in Theorem \ref{existence} with $
N\geq 5$. There exists a sufficiently small $\delta_0=\delta_0(N)$ such that if $\mathcal{ {E}}_N(0)\le \delta_0$,
and assuming that $(u_0,E_0,B_0)\in \dot{H}^{-s}$ for some $s\in [0,3/2)$ or $(u_{0},E_0,B_0)\in \dot{B}_{2,\infty}^{-s}$ for some $s\in (0,3/2]$,  then we have
\begin{equation}\label{H-sbound}
\norm{(u,E,B)(t)}_{\dot{H}^{-s}}\le C_0
\end{equation}or
\begin{equation}\label{H-sbound Besov}
\norm{(u,E,B)(t)}_{\dot{B}_{2,\infty}^{-s}}\le C_0.
\end{equation}
Moreover, for any fixed integer $k\ge 0$, if $N\ge 2k+2+s$, then
\begin{equation}\label{basic decay}
\norm{\na^k(n,u,E,B)(t)}_{L^2}\le  C_0 (1+ t)^{- \frac{k+s}{2} }.
\end{equation}

Furthermore, for any fixed integer $k\ge 0$, if $N\ge2k+4+s$, then
\begin{equation}\label{further decay1}
\norm{\na^k(n,u,E)(t)}_{L^2}\le  C_0 (1+ t)^{- \frac{k+1+s}{2} };
\end{equation}
if $N\ge2k+6+s$, then
\begin{equation}\label{further decay11}
\norm{\na^k n (t)}_{L^2}\le  C_0 (1+ t)^{- \frac{k+2+s}{2} };
\end{equation}
if $N\ge2k+12+s$ and $ B_\infty =0$, then
\begin{equation}\label{further decay2}
\norm{\na^k (n,{\rm div}u) (t)}_{L^2}\le  C_0 (1+ t)^{- (\frac k2+\frac74+s) }.
\end{equation}
\end{theorem}

The proof of Theorem \ref{decay} is based on the regularity
interpolation method developed in  Strain and Guo \cite{SG06}, Guo
and Wang \cite{GW} and Sohinger and Strain \cite{SS}. To prove the
optimal decay rate of the dissipative equations in the whole space,
Guo and Wang \cite{GW} developed a general energy method of using a
family of scaled energy estimates with minimum derivative counts and
interpolations among them. Note that the homogeneous Sobolev space
$\Dot{H}^{-s}$ was introduced there to enhance the decay rates. By
the usual embedding theorem, we know that for $p\in (1,2]$,
$L^p\subset \Dot{H}^{-s}$ with
$s=3(\frac{1}{p}-\frac{1}{2})\in[0,3/2)$. Hence the $L^p$--$L^2$
type of the optimal decay results follows as a corollary. However,
this does not cover the case $p=1$. To amend this, Sohinger and
Strain \cite{SS} instead introduced the homogeneous Besov space
$\Dot{B}_{2,\infty}^{-s}$ due to the fact that the endpoint
embedding  $L^1\subset \Dot{B}_{2,\infty}^{-\frac{3}{2}}$ holds. The
method of \cite{GW,SS} can be applied to many dissipative equations
in the whole space, however, it cannot be applied directly to the
compressible Euler-Maxwell system which is of regularity-loss. To
get around this difficulty, based on the refined energy estimates
stated in Lemma \ref{energy lemma}--\ref{other di}, we deduce
\begin{equation}
\frac{d}{dt} {\mathcal{E}}_k^{k+2}+\mathcal{D}_k^{k+2}\le C_k   \norm{(n,u)}_{L^\infty}\norm{\na^{k+2}(n, u )}_{L^2}
\norm{ \na^{k+2}( E,  B )}_{L^2},
\end{equation}
where ${\mathcal{E}}_k^{k+2}$ and $\mathcal{D}_k^{k+2}$ with minimum
derivative counts are defined by \eqref{1111} and \eqref{2222}
respectively. Then combining the methods of \cite{GW,SS} and a trick
of Strain and Guo \cite{SG06} to treat the electromagnetic field,
that is, doing the regularity interpolation in a double way, we are
able to conclude the decay rate \eqref{basic decay}. The decay rate
of $B$ in (1.16) is optimal as it is consistent with the linear one
proved in Duan [2]. Indeed, the decay rate of $B$ is the slowest
among all the components of the solution. In this sense, if in view
of the whole solution, we may regard (1.16) as to be optimal. The
faster decay rates \eqref{further decay1}--\eqref{further decay2}
follow by revisiting the equations carefully. In particular, we will
use a bootstrap argument to derive \eqref{further decay2}.

As quoted above, by Theorem \ref{decay}, we have the following corollary of the usual $L^p$--$L^2$ type of the decay results:
\begin{Corollary}\label{2mainth}
Under the assumptions of Theorem \ref{decay} except that we replace the $\dot{H}^{-s}$ or $\Dot{B}_{2,\infty}^{-s}$ assumption  by that $(u_{0},E_0,B_0)\in L^p$ for some $p\in [1,2]$, then for any fixed integer $k\ge 0$, if $N\ge 2k+2+ s_p $, then
\begin{equation}\label{p21}
\norm{\na^k(n,u,E,B)(t)}_{L^2}\le  C_0 (1+ t)^{- \frac{k+s_p}{2} }.
\end{equation}
Here the number $s_{p} :=3\left(\frac{1}{p}-\frac{1}{2}\right)$.

Furthermore, for any fixed integer $k\ge 0$, if $N\ge2k+4+s_p$, then
\begin{equation} \label{p22}
\norm{\na^k(n,u,E)(t)}_{L^2}\le  C_0 (1+ t)^{- \frac{k+1+s_p}{2} };
\end{equation}
if $N\ge2k+6+s_p$, then
\begin{equation} \label{p23}
\norm{\na^k n (t)}_{L^2}\le  C_0 (1+ t)^{- \frac{k+2+s_p}{2} };
\end{equation}
if $N\ge2k+12+s_p$ and $ B_\infty =0$, then
\begin{equation}\label{n L^1}
\norm{\na^k (n,{\rm div}u) (t)}_{L^2}\le  C_0 (1+ t)^{- (\frac k2+\frac74+s_p) }.
\end{equation}
\end{Corollary}

The following are several remarks for Theorem \ref{existence}, Theorem \ref{decay} and Corollary \ref{2mainth}.

\begin{Remark}
In Theorem \ref{existence}, we only assume the $H^3$ norm of the initial data is small, but the higher order derivatives can be arbitrarily large.
Notice that in Theorem \ref{decay} the $\dot{H}^{-s}$ and $\dot{B}_{2,\infty}^{-s}$ norms of the solution are preserved along the time evolution; however, in Corollary \ref{2mainth} it is difficult to show that the $L^p$ norm of the solution can be preserved.
Note that the $L^2$ decay rate of the higher order spatial derivatives of the solution
is obtained. Then the general optimal $L^q$ $(2\le q\le \infty)$ decay rates of the solution follow by the Sobolev
interpolation.
\end{Remark}

\begin{Remark}
We remark that Corollary \ref{2mainth} not only provides an alternative approach to derive the $L^p$--$L^2$ type
of the optimal decay results but also improves the previous results of the  $L^p$--$L^2$ approach in Ueda and Kawashima \cite{UK} and Duan \cite{D}.
In Ueda and Kawashima \cite{UK}, the decay rates \eqref{p21}--\eqref{p23} with $p=2$ were proved by using the time weighted
energy method, and when $p=1$ they were proved by combining the time weighted
energy method and the linear decay analysis but under the stronger assumption that  $\norm{(n_0,u_0,E_0,B_0)}_{L^1}$ is sufficiently small.
In Duan \cite{D}, assuming that $B_\infty=0$ and $\norm{(u_0,E_0,B_0)}_{L^1}$ is sufficiently small, by combining the energy method and the linear decay analysis, Duan proved that
\begin{equation}
\norm{n (t)}_{L^2}\le  C_0 (1+ t)^{- \frac{11}{4} },\quad \norm{(u,E)(t)}_{L^2}\le  C_0 (1+ t)^{- \frac{5}{4} } \text{ and } \norm{B (t)}_{L^2}\le  C_0 (1+ t)^{- \frac{3}{4} }.
\end{equation}
Note that we have removed the smallness of $L^p$ norm of initial
data, and for $p=1$ our decay rate of $n(t)$ in \eqref{n L^1} is
$(1+t)^{-13/4}$. Besides, Duan \cite{D} essentially depended on the
assumption $B_\infty=0$, and our results \eqref{p21}--\eqref{p23}
work for the general case $B_\infty\neq0$.
\end{Remark}

The rest of our paper is organized as follows. In section \ref{section2}, we establish the refined energy estimates for the solution and derive the negative Sobolev and Besov estimates. Theorem \ref{existence} and Theorem \ref{decay} are proved in section \ref{section3}.

%%%%%%%%%%%%%%%%%%%%%%%%%%%%%%%%%%%%%%%%%%%%%%%%%%%%%%%%%%%%%%%%%%%%%%%%%%%%%
\section{Nonlinear energy estimates}\label{section2}
%%%%%%%%%%%%%%%%%%%%%%%%%%%%%%%%%%%%%%%%%%%%%%%%%%%%%%%%%%%%%%%%%%%%%%%%%%%%%

In this section, we will do the a priori estimate by assuming that $\norm{n(t)}_{H^3}\le \delta\ll 1$. Recall the expression \eqref{f} of $f(n)$. Then by Taylor's formula and Sobolev's inequality, we have
\begin{equation}\label{fn}
f(n)\sim n\hbox{ and } \norms{f^{(k)}(n)} \le C_k\hbox{ for any }k\ge 1.
\end{equation}

%%%%%%%%%%%%%%%%%%%%%%%%%%%%%%%%%%%%%%%%%%%%%%%%%%%%%%%%%%%%%%%%%%%%%%%%%%%%%
\subsection{Preliminary}
%%%%%%%%%%%%%%%%%%%%%%%%%%%%%%%%%%%%%%%%%%%%%%%%%%%%%%%%%%%%%%%%%%%%%%%%%%%%%

In this subsection, we collect the analytic tools which will be used in the paper and prove a basic estimate for the nonlinear function $f(n)$.

 \begin{lemma}\label{A1}
Let $2\le p\le +\infty$ and $\alpha,m,\ell\ge 0$. Then we have
\begin{equation}\label{A.1}
\norm{\na^\alpha f}_{L^p}\le C_{p} \norm{ \na^mf}_{L^2}^{1-\theta}
\norm{ \na^\ell f}_{L^2}^{\theta}.
\end{equation}
Here $0\le \theta\le 1$ (if $p=+\infty$, then we require that $0<\theta<1$) and $\alpha$ satisfy
\begin{equation}\label{A.2}
\alpha+3\left(\frac12-\frac{1}{p}\right)=m(1-\theta)+\ell\theta.
\end{equation}
\end{lemma}
\begin{proof}
For the case $2\le p<+\infty$, we refer to Lemma A.1 in \cite{GW}; for the case $p=+\infty$, we refer to Exercise 6.1.2 in \cite{Gla} (pp. 421).
\end{proof}

\begin{lemma}\label{A2}
For any integer $k\ge0$, we have
\begin{equation}\label{fkn0}
\norm{\na^kf(n)}_{L^\infty} \le C_k\norm{\na^k n}
_{L^2}^{1/4}\norm{\na^{k+2}n}_{L^2} ^{3/4},
\end{equation}
and
\begin{equation}\label{fkn}
\norm{\na^kf(n)}_{L^2}\le C_k\norm{\na^kn}_{L^2}.
\end{equation}
\end{lemma}
\begin{proof}
The proof is based on Lemma \ref{A1}. For \eqref{fkn0}, we refer to Lemma 3.1 in \cite{GW}. For \eqref{fkn}, in light of \eqref{fn}, it suffices to prove that when $k\ge 1$, \eqref{fkn} holds for all $f(n)$ with bounded derivatives.
We will use an induction on $k\ge 1$. If $k=1$, we have
\begin{equation}
\norm{\na f(n)}_{L^2}=\norm{ f'(n)\na n}_{L^2}\lesssim \norm{\na n}_{L^2}.
\end{equation}
Assume \eqref{fkn} holds for from $1$ to $k-1$.  We use the Leibniz formula to have
\begin{equation}\label{fkn er}
\begin{split}
 \norm{\na^{k}f(n)}_{L^2}&=\norm{\na^{k-1}(f'(n)\na n)}_{L^2}
 \\&
 \le C_k\left(\norm{ f'(n)\na^{k} n}_{L^2}+\norm{\na f'(n)\na^{k-1} n}_{L^2}+\sum_{ \ell=2}^{ k-1}\norm{\na^{\ell}f'(n)\na^{k-\ell} n}_{L^2}\right)
 .
 \end{split}
 \end{equation}
Here if $k=2$, then the summing term in \eqref{fkn er} is nothing, etc. By H\"older's inequality and Sobolev's inequality, we have
\begin{equation}
\norm{ f'(n)\na^{k} n}_{L^2}+\norm{\na f'(n)\na^{k-1} n}_{L^2}
\lesssim \norm{ \na^{k} n}_{L^2}+\norm{\na n}_{L^3}\norm{\na^{k-1} n}_{L^6}\lesssim \norm{ \na^{k} n}_{L^2}.
\end{equation}
For the summing term we use the induction hypothesis to obtain that for $2\le \ell\le k-1$,
 \begin{equation}
\norm{\na^{\ell}f'(n)\na^{k-\ell} n}_{L^2} \le \norm{\na^{\ell}f'(n)}_{L^2}\norm{\na^{k-\ell} n}_{L^\infty}
\lesssim \norm{\na^{\ell} n}_{L^2}\norm{\na^{k-\ell} n}_{L^\infty}.
  \end{equation}
By Lemma \ref{A1}, if $\ell\le \left[\frac{k-1}{2}\right]$, then we have
   \begin{equation}
   \norm{\na^{\ell} n}_{L^2}\norm{\na^{k-\ell} n}_{L^\infty}
 \lesssim \norm{\na^\alpha n}_{L^2}^{\frac{k-\ell+\frac{3}{2}}{k}}\norm{\na^{k} n}_{L^2}^{\frac{\ell-\frac{3}{2}}{k}}
\norm{ n}_{L^2}^{\frac{\ell-\frac{3}{2}}{k}}\norm{\na^{k} n}_{L^2}^{\frac{k-\ell+\frac{3}{2}}{k}}
 \lesssim  \norm{\na^{k} n}_{L^2},
 \end{equation}
 where $\al$ is defined by
 \begin{equation}
\ell=\al\times \frac{k-\ell+\frac{3}{2}}{k} +k\times\frac{\ell-\frac{3}{2}}{k}\Longrightarrow \al=\frac{3k}{2(k-\ell)+3}<3;
\end{equation}
if $\ell\ge \left[\frac{k-1}{2}\right]+1$, then we have
   \begin{equation}
   \norm{\na^{\ell} n}_{L^2}\norm{\na^{k-\ell} n}_{L^\infty}
 \lesssim \norm{ n}_{L^2}^{1-\frac{\ell}{k}}\norm{\na^{k} n}_{L^2}^\frac{\ell}{k}
\norm{\na^\alpha n}_{L^2}^{ \frac{\ell}{k}}\norm{\na^{k} n}_{L^2}^{1-\frac{\ell}{k}}
 \lesssim  \norm{\na^{k} n}_{L^2},
 \end{equation}
 where $\al$ is defined by
 \begin{equation}
k-\ell+ \frac{3}{2}=\al\times \frac{\ell}{k} +k\times\left(1-\frac{\ell}{k}\right)\Longrightarrow \al=\frac{3k}{2\ell}<3.
\end{equation}
We thus conclude the lemma.
\end{proof}

We recall the following commutator estimate:
\begin{lemma}\label{commutator}
Let $k\ge 1$ be an integer and define the commutator
\begin{equation}\label{commuta}
\left[\na^k,g\right]h=\na^k(gh)-g\na^kh.
\end{equation}
Then we have
\begin{equation}
\norm{\left[\na^k,g\right]h}_{L^2} \le C_k\left( \norm{\na g}_{L^\infty}%
\norm{\na^{k-1}h}_{L^2}+\norm{\na^k g}_{L^2}\norm{ h}_{L^\infty}\right).
\end{equation}
\end{lemma}
\begin{proof}
It can be proved by using Lemma \ref{A1}, see Lemma 3.4 in \cite{MB} (pp. 98) for instance.
\end{proof}

We have the $L^p$ embeddings:
\begin{lemma}\label{Riesz lemma}
Let $0\le s<3/2,\ 1<p\le 2$ with $1/2+s/3=1/p$, then
\begin{equation}\label{Riesz estimate}
\norm{ f}_{\dot{H}^{-s}}\lesssim\norm{ f}_{L^p}.
\end{equation}
\end{lemma}
\begin{proof}
It follows from the Hardy-Littlewood-Sobolev theorem, see \cite{Gla}.
\end{proof}

\begin{lemma}\label{Lp embedding}
Let $0< s\le 3/2,\ 1\le p<2$ with $1/2+s/3=1/p$, then
\begin{equation}\label{Lp embedding inequality}
\norm{f}_{\dot{B}_{2,\infty}^{-s}}\lesssim\norm{f}_{L^p}.
\end{equation}
\end{lemma}
\begin{proof}
See Lemma 4.6 in \cite{SS}.
\end{proof}

It is important to use the following special interpolation estimates:
\begin{lemma}\label{1-sinte}
Let $s\ge 0$ and $\ell\ge 0$, then we have
\begin{equation}  \label{1-sinterpolation}
\norm{\na^\ell f}_{L^2}\le \norm{\na^{\ell+1} f}_{L^2}^{1-\theta}%
\norm{ f}_{\dot{H}^{-s}}^\theta, \hbox{ where }\theta=\frac{1}{\ell+1+s}.
\end{equation}
\end{lemma}
\begin{proof}
It follows directly by the Parseval theorem  and H\"older's
inequality.
\end{proof}

\begin{lemma}\label{Besov interpolation}
Let $s> 0$ and $\ell\ge 0$, then we have
\begin{equation}
\norm{\na^\ell f}_{L^2}\le \norm{\na^{\ell+1} f}_{L^2}^{1-\theta}%
\norm{ f}_{\dot{B}^{-s}_{2,\infty}}^\theta, \hbox{ where }\theta=\frac{1}{\ell+1+s}.
\end{equation}
\end{lemma}
\begin{proof}
See Lemma 4.5 in \cite{SS}.
\end{proof}

%%%%%%%%%%%%%%%%%%%%%%%%%%%%%%%%%%%%%%%%%%%%%%%%%%%%%%%%%%%%%%%%%%%%%%%%%%%%%
\subsection{Energy estimates}
%%%%%%%%%%%%%%%%%%%%%%%%%%%%%%%%%%%%%%%%%%%%%%%%%%%%%%%%%%%%%%%%%%%%%%%%%%%%%

In this subsection, we will derive the basic energy estimates for the solution to the Euler-Maxwell system \eqref{EM per}. We begin with the standard energy estimates.
\begin{lemma}\label{energy lemma}
For any integer $k\ge0$, we have
\begin{align}\label{energy 1}
&\frac{d}{dt}\sum_{l=k}^{k+2}\norm{\na^{l} ( n, u, E, B )}_{L^2}^2 +\la\sum_{l=k}^{k+2}\norm{\na^{l} u}_{L^2}^2 \nonumber\\
&\quad\le C_k  \left(\norm{ (n, u)}_{ H^{k+1}\cap H^{\frac{k}{2}+2}\cap H^3}+\norm{\na B}_{L^2}  \right)
\left(\sum_{l=k}^{k+2}\norm{\na^{l} (n, u )}_{L^2}^2+\sum_{l=k}^{k+1}\norm{\na^{l}E}_{L^2}^2+\norm{\na^{k+1} B}_{L^2}^2\right)
\nonumber\\&\qquad+\norm{(n,u)}_{L^\infty}\norm{\na^{k+2} (n, u )}_{L^2} \norm{\na^{k+2} ( E,  B )}_{L^2}.
\end{align}
\end{lemma}
\begin{proof}
The standard $\na^l$ ($l=k,k+1,k+2$) energy estimates on the system $\eqref{EM per}$ yield
\begin{equation}\label{I1-I4}
\begin{split}
&\frac{1}{2}\frac{d}{dt}\int \norms{\na^{l}(n,u,E,B)}^2  +\nu  \norm{\na^l u}_{L^2}^2 \\
&\quad=-\mu \int \na^l (n{\rm div} u)\na^ln+\na^l(n\na n)\cdot\na^l u
-\int \na^l (u\cdot\na n)\na^ln+ \na^l(u\cdot\na u)\cdot\na^l u \\
&\qquad-\int \na^l (u\times B)\cdot\na^l u +\nu \int \na^{l} (f(n)u)\cdot\na^{l} E \\
&\quad:=I_1+I_2+I_3+ \nu I_4.
\end{split}
\end{equation}

We now estimate $I_1- I_4$. First, we use the commutator notation \eqref{commuta} to rewrite $I_1$ as
\begin{equation}
\begin{split}
I_1&=-\mu \int   \left( n{\rm div} \na^l u+ \left[\na^l ,n \right]{\rm div}  u\right)\na^ln
+ \left(n\na \na^l n +\left[\na^l ,n \right]\na  n\right) \cdot \na^l u
\\& =-\mu \int   n{\rm div}( \na^l u \na^l n)+ \left[\na^l ,n \right]{\rm div}  u \na^ln
 +\left[\na^l ,n \right]\na  n  \cdot \na^l u.
\end{split}
\end{equation}
By integrating by parts, we have
\begin{equation}
 -\int   n{\rm div}( \na^l u \na^l n)  =\int  \na n\cdot \na^l u \na^l n
\le \norm{\na n}_{L^\infty}\norm{\na^l u}_{L^2}\norm{\na^l n}_{L^2}.
\end{equation}
We use the commutator estimate of Lemma \ref{commutator} to bound
\begin{equation}
\begin{split}
-\int \left [\na^l ,n \right]{\rm div}  u \na^ln
 &\le C_l\left( \norm{\na n}_{L^\infty}
\norm{\na^{l-1}{\rm div}  u }_{L^2}+\norm{\na^l n}_{L^2}\norm{{\rm div}  u}_{L^\infty}\right)\norm{\na^l n}_{L^2}
\\& \le C_l\norm{\na (n,u)}_{L^\infty}\norm{\na^l (n,u)}_{L^2}\norm{\na^l n}_{L^2}.
\end{split}
\end{equation}
Bounding the last term of $I_1$ similarly, and then applying the same arguments to the term $I_2$, by Sobolev's  and Cauchy's inequalities, we deduce
\begin{equation}\label{I1 I2}
I_1+I_2\le  C_l \norm{ ( n,u)}_{H^3}\norm{\na^l (n,u)}_{L^2}^2.
\end{equation}

Next, we estimate the term $I_3$, and we must be much more careful with this term since the magnetic field $B$ has the weakest dissipative estimates. First of all, we have
\begin{equation}\label{B yi}
\begin{split}
I_{3}&=-\sum_{\ell=1}^{ l} C_l^\ell \int \na^{l-\ell} u\times \na^{\ell} B \cdot \na^l u\le C_l\sum_{\ell=1}^{l}\norm{\na^{l-\ell}u \na^{\ell}B}_{L^2}\norm{\na^l u}_{L^2}.
\end{split}
\end{equation}
Here, if $l<1$, then it's nothing, and etc.
We have to distinguish the arguments by the value of $l$. First, let $l=k$. We take $L^3-L^6$ and then apply Lemma \ref{A1} to have
\begin{equation}
\begin{split}
\norm{\na^{k-\ell} u\na^{\ell}B}_{L^{2}}&\lesssim \norm{\na^{k-\ell} u}_{L^3} \norm{\na^{\ell}B}_{L^6}
\\&\lesssim \norm{\na^\alpha u}^{\frac {\ell}{k}}_{L^2}\norm{\na^{k}u}_{L^2}^{1-\frac{\ell}{k}} \norm{\na B}_{L^2}^{1-\frac{\ell}{k}}\norm{\na^{k+1}B}_{L^2}^{\frac {\ell}{k}},
 \end{split}
\end{equation}
where $\alpha$ is defined by
\begin{equation}
k-\ell+\frac12=\alpha\times\frac{\ell}{k}
+k\times \left(1-\frac{\ell}{k}\right)
 \Longrightarrow \alpha=\frac{ k }{2\ell}\le \frac{k}{2}.
\end{equation}
Hence by Young's inequality, we have that for $l=k$,
\begin{equation}\label{I3 k}
I_3\le C_k\left( \norm{ u}_{H^{\frac{k}{2}}}  +\norm{\na B}_{L^2}  \right) \left(\norm{\na^{k}u}_{L^2}^2 +\norm{\na^{k+1}B}_{L^2}^2 \right).
\end{equation}

We then let $l=k+1$. If $1\le \ell\le k$, we take $L^3-L^6$ and by  Lemma \ref{A1} again to obtain
\begin{equation}
\begin{split}
\norm{\na^{k+1-\ell} u\na^{\ell}B}_{L^{2}}&\lesssim \norm{\na^{k+1-\ell} u}_{L^3} \norm{\na^{\ell}B}_{L^6}
\\&\lesssim \norm{\na^\alpha u}^{\frac {\ell}{k}}_{L^2}\norm{\na^{k+1}u}_{L^2}^{1-\frac{\ell}{k}} \norm{\na B}_{L^2}^{1-\frac{\ell}{k}}\norm{\na^{k+1}B}_{L^2}^{\frac {\ell}{k}},
\end{split}
\end{equation}
where $\alpha$ is defined by
\begin{equation}
k+1-\ell+\frac12=\alpha\times\frac {\ell}{k}
+(k+1)\times \left({1-\frac{\ell}{k}}\right)
\Longrightarrow \alpha=1+\frac{k}{2\ell}\le \frac{k}{2}+1;
\end{equation}
if $\ell=k+1$, we take $L^\infty-L^2$ to get
\begin{equation}
\norm{ u\na^{k+1}B}_{L^{2}}\lesssim \norm{ u}_{L^\infty} \norm{\na^{k+1}B}_{L^2}.
\end{equation}
We thus have that for $l=k+1$, by Sobolev's inequality,
\begin{equation}
I_3\le C_k\left( \norm{  u}_{H^{\frac{k}{2}+1}\cap H^2}  +\norm{\na B}_{L^2}  \right) \left(\norm{\na^{k+1}u}_{L^2}^2 +\norm{\na^{k+1}B}_{L^2}^2 \right).
\end{equation}

We now let $l=k+2$. If $1\le \ell\le k$, we take $L^3-L^6$ and by Lemma \ref{A1} again to have
\begin{equation}
\begin{split}
\norm{\na^{k+2-\ell} u\na^{\ell}B}_{L^{2}}&\lesssim \norm{\na^{k+2-\ell} u}_{L^3} \norm{\na^{\ell}B}_{L^6}
\\&\lesssim \norm{\na^\alpha u}^{\frac {\ell}{k}}_{L^2}\norm{\na^{k+2}u}_{L^2}^{1-\frac{\ell}{k}} \norm{\na B}_{L^2}^{1-\frac{\ell}{k}}\norm{\na^{k+1}B}_{L^2}^{\frac {\ell}{k}},
\end{split}
\end{equation}
where $\alpha$ is defined by
\begin{equation}
k+2-\ell+\frac12=\alpha\times\frac {\ell}{k}
+(k+2)\times \left({1-\frac{\ell}{k}}\right)
\Longrightarrow \alpha=2+\frac{k}{2\ell}\le \frac{k}{2}+2;
\end{equation}
if $\ell=k+1$ or $ k+2$, we take $L^\infty-L^2$ to get
\begin{equation}
\norm{ \na u\na^{k+1}B}_{L^{2}}\lesssim \norm{\na u}_{L^\infty} \norm{\na^{k+1}B}_{L^2},
\end{equation}
and
\begin{equation}
\norm{ u\na^{k+2}B}_{L^{2}}\lesssim \norm{  u}_{L^\infty} \norm{\na^{k+2}B}_{L^2}.
\end{equation}
We thus have that for $l=k+2$,
\begin{equation}
\begin{split}
I_3&\le C_k\left( \norm{  u}_{H^{\frac{k}{2}+2}\cap H^3}  +\norm{\na B}_{L^2}  \right) \left(\norm{\na^{k+2}u}_{L^2}^2 +\norm{\na^{k+1}B}_{L^2}^2 \right)
\\&\quad+C\norm{ u}_{L^\infty} \norm{\na^{k+2}B}_{L^2}\norm{\na^{k+2}u}_{L^2}.
\end{split}
\end{equation}

We now estimate the last term $I_4$. We again have to distinguish the arguments by the value of $l$. First, for $l=k$ or $k+1$, we have
\begin{equation}
I_{4}=\sum_{\ell=0}^{ l}C_l^\ell \int \na^{\ell}f(n)\na^{l-\ell}u\cdot\na^{l} E \le C_l \sum_{\ell=0}^{ l}\norm{\na^{\ell}f(n)\na^{l-\ell}u}_{L^2}\norm{\na^{l} E}_{L^2}.
\end{equation}
If $0\le \ell\le l-2$, we take $L^\infty-L^2$ and by Lemma \ref{A1} and the estimate \eqref{fkn0} of Lemma \ref{A2} to obtain
\begin{equation}
\begin{split}
\norm{\na^{\ell}f(n)\na^{l-\ell}u}_{L^2}&\le  \norm{\na^{\ell} f(n)}_{L^\infty} \norm{\na^{l-\ell}u}_{L^2}
\\&\le C_l\norm{\na^{\ell} n}_{L^2}^{\frac14} \norm{\na^{\ell+2}n}_{L^2}^{\frac34}\norm{\na^{l-\ell}u}_{L^2}
\\&\le C_l \left(\norm{n}_{L^2}^{1-\frac{\ell}{l}}\norm{\na^{l} n}_{L^2}^{\frac{\ell}{l}}\right)^{\frac14}  \left(\norm{ n}_{L^2}^{1-\frac{\ell+2}{l}} \norm{\na^{l} n}_{L^2}^{\frac{\ell+2}{l}}\right)^{\frac34}\norm{\na^{l-\ell}u}_{L^2}
\\&\le C_l \norm{n}_{L^2}^{1-\frac{2\ell+3}{2l}}\norm{\na^{l} n}_{L^2}^{\frac{2\ell+3}{2l}} \norm{\na^{\alpha} u}_{L^2}^{\frac{2\ell+3}{2l}} \norm{\na^{l} u}_{L^2}^{1-\frac{2\ell+3}{2l}},
\end{split}
\end{equation}
where $\alpha$ is defined by
\begin{equation}l-\ell=\alpha\times\frac{2\ell+3}{2l}
+l\times\left(1-\frac{2\ell+3}{2l}\right)
\Longrightarrow \alpha=\frac{3l}{2\ell+3}\le l;
\end{equation}
if $\ell=l-1$, we take $L^6-L^3$ and by the estimate \eqref{fkn} of Lemma \ref{A2} to have
\begin{equation}
\norm{\na^{l-1}f(n)\na u}_{L^2}\le  \norm{\na^{l-1}f(n)}_{L^6}\norm{\na u}_{L^3}  \le C_l \norm{\na^{l} n }_{L^2}\norm{  u}_{H^2};
\end{equation}
if $\ell=l$, we take $L^2-L^\infty$ and by the estimate \eqref{fkn0} of Lemma \ref{A2} to have
\begin{equation}
\norm{\na^{l}f(n)u}_{L^2}\le \norm{\na^{l}f(n)}_{L^2}\norm{u}_{L^\infty}
\le C_l\norm{\na^{l} n }_{L^2}\norm{u}_{H^2}.
\end{equation}
We thus have that for $l=k$ or $k+1$,
\begin{equation}\label{I4 k k+1}
I_4\le C_l  \norm{  (n,u)}_{H^{l}\cap H^2}  \left(\norm{\na^l (n, u)}_{L^2}^2 +\norm{\na^l E}_{L^2}^2 \right).
\end{equation}

Now for $l=k+2$, we rewrite $I_4$ as
\begin{equation}
\begin{split}
I_{4}&=\sum_{\ell=0}^{ k+2}C_{k+2}^{\ell}\int \na^{\ell}f(n)\na^{k+2-\ell}u\cdot\na^{k+2} E \\
&=\int \left(f(n)\na^{k+2}u+ \na^{k+2}f(n)u\right)\cdot\na^{k+2} E
\\
&\quad-\sum_{\ell=1}^{ k+1} C_{k+2}^{\ell}\int \na\left(\na^{k+2-\ell}f(n)\na^{\ell}u\right)\cdot\na^{k+1} E
\\
&=\int \left(f(n)\na^{k+2}u+ \na^{k+2}f(n)u\right)\cdot\na^{k+2} E\\
&\quad-(k+2)\int\left(\na^{k+2}f(n)\na u+\na f(n)\na^{k+2}u\right)\cdot\na^{k+1} E
\\&\quad-\sum_{\ell=2}^{ k+1} C_{k+2}^{\ell}\int \na^{k+3-\ell}f(n)\na^{\ell}u \cdot\na^{k+1} E
-\sum_{\ell=1}^{ k} C_{k+2}^{\ell}\int  \na^{k+2-\ell}f(n)\na^{\ell+1}u\cdot\na^{k+1} E
\\&:=I_{41} +I_{42} +I_{43} .
\end{split}
\end{equation}
By Lemma \ref{A2}, we have
\begin{equation}
\begin{split}
I_{41} &\le C_k\left(\norm{f(n)}_{L^\infty}\norm{\na^{k+2}u}_{L^2}+\norm{\na^{k+2}f(n)}_{L^2}\norm{u}_{L^\infty}\right)\norm{\na^{k+2} E}_{L^2}
\\&\le C_k\norm{(n,u)}_{L^\infty}\norm{\na^{k+2}(n,u)}_{L^2}\norm{\na^{k+2} E}_{L^2}
\end{split}
\end{equation}
and
\begin{equation}
\begin{split}
I_{42} &\le C_k\left(\norm{\na^{k+2} f(n)}_{L^2}\norm{\na u}_{L^\infty}+\norm{\na f(n)}_{L^\infty}\norm{\na^{k+2} u}_{L^2}\right)\norm{\na^{k+1} E}_{L^2}
\\&\le C_k\norm{\na(n,u)}_{L^\infty}\norm{\na^{k+2}(n,u)}_{L^2}\norm{\na^{k+1} E}_{L^2}.
\end{split}
\end{equation}
As for the cases $l=k,k+1$ for $I_4$, we can bound $I_{43}$ by
\begin{equation}
\begin{split}
 I_{43} \le C_k\norm{(n,u)}_{H^{k+1}\cap H^2}  \left(\norm{\na^{k+1}( n,u)}_{L^2}^2 +\norm{\na^{k+1} E}_{L^2}^2 \right).
\end{split}
\end{equation}
Hence, we have that for $l=k+2$,
\begin{equation}
\begin{split}
 I_{4} &\le C_k\norm{(n,u)}_{H^{k+1}\cap H^3}  \left(\norm{\na^{k+1}( n,u)}_{L^2}^2 +\norm{\na^{k+2}( n,u)}_{L^2}^2+\norm{\na^{k+1} E}_{L^2}^2 \right)
 \\&\quad+C_k\norm{(n,u)}_{L^\infty}\norm{\na^{k+2}(n,u)}_{L^2}\norm{\na^{k+2} E}_{L^2}.
 \end{split}
\end{equation}

Consequently, plugging the estimates for $I_1- I_4$ into \eqref{I1-I4} with $l=k,k+1,k+2$, and then summing up, we deduce \eqref{energy 1}.
\end{proof}

Note that in Lemma \ref{energy lemma} we only derive the dissipative estimate of $u$. We now recover the dissipative
estimates of $ n, E$ and $B$ by constructing some interactive energy functionals in the following lemma.

\begin{lemma}\label{other di}
For any integer $k\ge0$, we have that for any small fixed $\eta>0$,
\begin{equation}\label{other dissipation}
\begin{split}
&\frac{d}{dt}\left(\sum_{l=k}^{k+1}\int \na^lu\cdot\na\na^{l} n +\sum_{l=k}^{k+1}\int \na^{l}u\cdot\na^lE -\eta\int
\na^k E \cdot\na^{k}\na\times B \right)
\\&\quad+\la \left(\sum_{l=k}^{k+2}\norm{\na^{l}n}_{L^2}^2+\sum_{l=k}^{k+1}\norm{\na^{l}E}_{L^2}^2+ \norm{\na^{k+1} B}_{L^2}^2\right)
\\ &\qquad\le
 C\sum_{l=k}^{k+2}\norm{\na^{l} u}_{L^2}^2
 +C_k  \left(\norm{ (n, u)}_{H^{k+1}\cap H^3}^2+\norm{\na B}_{L^2}^2  \right)  \left(\sum_{l=k}^{k+2}\norm{ \na^{l}(n, u )}_{L^2}^2+\norm{\na^{k+1}B}_{L^2}^2 \right).
\end{split}
\end{equation}
\end{lemma}

\begin{proof} We divide the proof into several steps.

{\it Step 1: Dissipative estimate of $n$.}

Applying $\na^l$ ($l=k,k+1$) to $\eqref{EM per}_2$ and then taking the $L^2$ inner product with $\na\na^{l} n$, we obtain
\begin{equation}  \label{lemmasanyi}
\begin{split}
\int  \partial_t\na^lu \cdot\na \na^{l} n +\norm{\na\na^{l } n}_{L^2}^2
&\le -\nu \int  \na^{l}E\cdot\na \na^{l} n +C \norm{\na^{l} u}_{L^2}\norm{\na^{l+1}  n}_{L^2}\\
&\quad  +\norm{\na^{l}\left(u\cdot\na
u+\mu n\na n+u\times B\right)}_{L^2}\norm{\na^{l+1}  n}_{L^2}
.
\end{split}
\end{equation}

The delicate first term on the left-hand side of \eqref{lemmasanyi}
involves $\partial_t\na^{l} u$, and the key idea is to integrate by parts
in the $t$-variable and use the continuity equation $\eqref{EM per}_1$. Thus integrating by
parts for both the $t$- and $x$-variables, we obtain
\begin{equation}\label{nt}
\begin{split}
\int  \na^{l} \partial_tu\cdot\na\na^l n&=\frac{d}{dt}\int
\na^{l}u\cdot\na\na^l n -\int  \na^{l}
 u\cdot \na \na^{l}\partial_t n
 \\
&=\frac{d}{dt}\int
\na^{l}u\cdot\na\na^l n +\int  \na^{l}
{\rm div} u\na^{l}\partial_t n  \\
& =\frac{d}{dt}\int
\na^{l}u\cdot\na\na^l n -\norm{\na^{l} {{\rm div} }u}_{L^2}^2-\int \na^{l} {\rm div} u\na^{l}\left(
u \cdot \na n+ \mu n{\rm div} u\right)
\\
& \ge \frac{d}{dt}\int
\na^{l}u\cdot \na\na^l n -C\norm{\na^{l+1} u}_{L^2}^2
 -C\norm{\na^{l}(u \cdot \na n)}_{L^2}^2-C\norm{\na^{l}( n {\rm div} u)}_{L^2}^2.
\end{split}
\end{equation}
Using the commutator estimate of Lemma \ref{commutator}, we have
\begin{equation}\label{e1}
\begin{split}
\norm{\na^{l}(u \cdot \na n)}_{L^2}&\le \norm{ u \cdot \na^{l} \na n }_{L^2}+\norm{\left[\na^l ,u \right]\cdot \na n }_{L^2}
\\&\le  \norm{ u}_{L^\infty} \norm{ \na^{l+1}   n }_{L^2}+C_l\norm{\na u}_{L^\infty}
\norm{\na^{l}n }_{L^2}+C_l\norm{\na^l u}_{L^2}\norm{\na n}_{L^\infty}
\\&\le C_l \norm{(n, u)}_{H^3} \left(\norm{ \na^{l}  ( n,u) }_{L^2}+\norm{ \na^{l+1}   n }_{L^2}\right).
\end{split}
\end{equation}
Similarly,
\begin{equation}\label{e2}
\norm{\na^{l}( n {\rm div} u)}_{L^2}\le C_l \norm{(n, u)}_{H^3} \left(\norm{ \na^{l}  ( n,u) }_{L^2}+\norm{ \na^{l+1}  u }_{L^2}\right).
\end{equation}
Hence, we obtain
\begin{equation}  \label{inter 1}
\begin{split}
 \int  \na^{l} \partial_tu\cdot\na\na^l n
&\ge \frac{d}{dt}\int
\na^{l}u\cdot\na^l\na n -C\norm{\na^{l+1} u}_{L^2}^2
\\ &\quad-C_l \norm{(n, u)}_{H^3}^2  \left(\norm{ \na^{l}  ( n,u) }_{L^2}^2+\norm{ \na^{l+1}  (n,u) }_{L^2}^2\right).
\end{split}
\end{equation}

Next, integrating by parts and using the equation $\eqref{EM per}_5$, we have
\begin{equation}\label{inter 3}
\begin{split}
-\nu \int  \na^{l}E\cdot \na\na^{l} n &=\nu \int  \na^{l}{\rm div} E\na^{l} n
=-\nu ^2\int  \na^{l}f(n)\na^{l} n \\&=-\nu ^2\int  \na^{l}(n+f(n)-n)\na^{l} n \\
&\le -\nu ^2\norm{\na^{l}n}_{L^2}^2+ C_l\norm{n}_{H^3}\norm{\na^{l}n}_{L^2}^2.
\end{split}
\end{equation}
Here we have used the estimate $\norm{\na^{l}(f(n)-n)}_{L^2}\le  C_l\norm{n}_{H^3}\norm{\na^{l}n}_{L^2}$, which follows by noticing that $f(n)-n\sim n^2$ and the similar arguments presented in Lemma \ref{A2}.

Lastly, as in \eqref{e1}--\eqref{e2}, we have
\begin{equation}
\begin{split}
\norm{\na^{l}\left(u\cdot\na
u+\mu n\na n \right)}_{L^2}\le C_l \norm{(n, u)}_{H^3}  \left(\norm{ \na^{l}  ( n,u) }_{L^2}+\norm{ \na^{l+1}  (n,u) }_{L^2}\right).
\end{split}
\end{equation}
From the estimate of $I_3$ in Lemma \ref{energy lemma}, we have that for $l=k$ or $k+1$,
\begin{equation} \label{inter 4}
\norm{\na^{l}\left(u\times B\right)}_{L^2} \le C_k\left( \norm{  u}_{H^{\frac{k}{2}+1}\cap H^2}  +\norm{\na B}_{L^2}  \right) \left(\norm{\na^l u}_{L^2} +\norm{\na^{k+1}B}_{L^2} \right).
\end{equation}

Plugging the estimates \eqref{inter 1}--\eqref{inter 4} into
\eqref{lemmasanyi}, by Cauchy's inequality, we obtain
\begin{equation}  \label{n estimate}
\begin{split}
\frac{d}{dt}&\sum_{l=k}^{k+1}\int \na^lu\cdot\na\na^{l} n +\la \sum_{l=k}^{k+2}\norm{\na^{l}n}_{L^2}^2
\\ &\le
 C\sum_{l=k}^{k+2}\norm{\na^{l} u}_{L^2}^2
 +C_k  \left(\norm{ (n, u)}_{H^{\frac{k}{2}+1}\cap H^3}^2+\norm{\na B}_{L^2}^2  \right)  \left(\sum_{l=k}^{k+2}\norm{ \na^{l}(n, u )}_{L^2}^2+\norm{\na^{k+1}B}_{L^2}^2 \right).
\end{split}
\end{equation}
This completes the dissipative estimate for $n$.

{\it Step 2: Dissipative estimate of $E$.}

Applying $\na^l$ ($l=k,k+1$) to $\eqref{EM per}_2$ and then taking the $L^2$ inner product with $\na^{l} E$, we obtain
\begin{equation}  \label{identy 1}
\begin{split}
\int  \na^l \partial_tu \cdot\na^{l}E +\nu \norm{\na^{l} E}_{L^2}^2
&\le -\int \na\na^{l} n\cdot \na^{l}E +C \norm{\na^{l} u}_{L^2}\norm{\na^{l}  E}_{L^2}\\
 &\quad +\norm{\na^{l}\left(u\cdot\na
u+\mu n\na n+u\times B\right)}_{L^2}\norm{\na^{l } E}_{L^2}.
\end{split}
\end{equation}

Again, the delicate first term on the left-hand side of \eqref{identy 1}
involves $\partial_t\na^{l} u$, and the key idea is to integrate by parts
in the $t$-variable and use the equation $\eqref{EM per}_3$ in the Maxwell system. Thus we obtain
\begin{equation}  \label{inter 21}
\begin{split}
 \int  \na^{l} \partial_tu\cdot\na^lE  &=\frac{d}{dt}\int
\na^{l}u\cdot\na^lE -\int  \na^{l}
u\cdot\na^{l}\partial_t E  \\
&=\frac{d}{dt}\int
\na^{l}u\cdot\na^lE -\nu \norm{\na^{l}u}_{L^2}^2-\nu \int \na^{l} u\cdot\na^{l}\left(f(n)u+
\na\times B\right) .
\end{split}
\end{equation}
From the estimates of $I_4$ in Lemma \ref{energy lemma}, we have that
\begin{equation}
\norm{\na^{l}\left(f(n)u\right)}_{L^2}\le C_l  \norm{  u}_{H^{l}\cap H^2}\norm{\na^l (n,u)}_{L^2} .
\end{equation}
We must be much more careful with the remaining term in \eqref{inter 21} since there is no small factor in front of it. The key is to use Cauchy's inequality and distinguish the cases of $l=k$ and $l=k+1$ due to the weakest dissipative estimate of $B$. For $l=k$, we have
\begin{equation}
-\nu \int \na^{k} u\cdot
\na\times \na^{k} B   \le \varepsilon \norm{ \na^{k+1} B}_{L^2}^2+C_\varepsilon \norm{\na^{k} u}_{L^2}^2;
\end{equation}
for $l=k+1$, integrating by parts, we obtain
\begin{equation}\label{inter 24}
\begin{split}
-\nu \int \na^{k+1} u\cdot
\na\times \na^{k+1} B  =
-\nu \int \na\times \na^{k+1} u\cdot
 \na^{k+1} B
\le \varepsilon \norm{ \na^{k+1} B}_{L^2}^2+C_\varepsilon \norm{\na^{k+2} u}_{L^2}^2.
\end{split}
\end{equation}

Plugging the estimates \eqref{inter 21}--\eqref{inter 24} and \eqref{inter 3}--\eqref{inter 4} from Step 1 into
\eqref{identy 1}, by Cauchy's inequality, we then obtain
\begin{equation}  \label{E estimate}
\begin{split}
& \frac{d}{dt}\sum_{l=k}^{k+1}\int \na^{l}u\cdot\na^lE +\la \sum_{l=k}^{k+1}\norm{\na^{l}E}_{L^2}^2
\\ &\quad \le \varepsilon \norm{ \na^{k+1} B}_{L^2}^2+
 C_\varepsilon \sum_{l=k}^{k+2}\norm{\na^{l} u}_{L^2}^2
  \\ &\qquad +C_k  \left(\norm{ (n, u)}_{H^{k+1}\cap H^3}^2+\norm{\na B}_{L^2}^2  \right)  \left(\sum_{l=k}^{k+2}\norm{ \na^{l}(n, u )}_{L^2}^2+\norm{\na^{k+1}B}_{L^2}^2 \right).
\end{split}
\end{equation}
This completes the dissipative estimate for $E$.

{\it Step 3: Dissipative estimate of $B$.}

Applying $\na^k$ to $\eqref{EM per}_3$ and then taking the $L^2$ inner product with $-\na\times\na^{k} B$,  we obtain
\begin{equation}  \label{D12}
\begin{split}
&-\int \na^k \partial_tE \cdot\na\times\na^{k} B +\nu \norm{\na\times\na^{k} B}_{L^2}^2
\\&\quad\le \nu \norm{\na^{k}u}_{L^2}\norm{\na\times\na^{k} B}_{L^2}+
 \nu \norm{\na^{k}(f(n)u)}_{L^2}\norm{\na\times\na^{k} B}_{L^2}.
\end{split}
\end{equation}

Integrating by parts for both the $t$- and $x$-variables and using the equation $\eqref{EM per}_4$, we have
\begin{equation}  \label{inter 31}
\begin{split}
-\int \na^k \partial_tE \cdot\na\times\na^{k} B
&=-\frac{d}{dt}\int
\na^k E \cdot\na\times\na^{k} B +\int  \na\times\na^{k}
E\cdot\na^{k}\partial_t B  \\
&=-\frac{d}{dt}\int
\na^k E \cdot\na\times\na^{k} B -\nu \norm{\na\times\na^{k}
E}_{L^2}^2.
\end{split}
\end{equation}
From the estimates of $I_4$ in Lemma \ref{energy lemma}, we have that
\begin{equation}\label{inter 32}
\norm{\na^{k}\left(f(n)u\right)}_{L^2}\le C_k  \norm{  u}_{H^{k}\cap H^2} \norm{\na^k (n,u)}_{L^2}.
\end{equation}

Plugging the estimates \eqref{inter 31}--\eqref{inter 32} into \eqref{D12} and by Cauchy's inequality, since ${\rm div} B=0$, we
then obtain
\begin{equation}   \label{B estimate}
\begin{split}
&-\frac{d}{dt}\int
\na^k E \cdot\na^{k}\na\times B +\la \norm{\na^{k+1} B}_{L^2}^2
\\&\quad\le C\norm{\na^{k}u}_{L^2}^2+C\norm{\na^{k+1}
  E}_{L^2}^2+
C_k  \norm{  u}_{H^{k}\cap H^2}^2 \norm{\na^k (n,u)}_{L^2}^2.
\end{split}
\end{equation}
This completes the dissipative estimate for $B$.

{\it Step 4: Conclusion.}

Multiplying \eqref{B estimate} by a small enough but fixed constant $\eta$ and then adding it with \eqref{E estimate}  so that the second term on the right-hand side of \eqref{B estimate} can be absorbed, then choosing $\varepsilon$ small enough so that the first term in \eqref{E estimate} can be absorbed; we obtain
\begin{equation}
\begin{split}
& \frac{d}{dt}\left(\sum_{l=k}^{k+1}\int \na^{l}u\cdot\na^lE-\eta\int
\na^k E \cdot\na^{k}\na\times B \right)+\la \left(\sum_{l=k}^{k+1}\norm{\na^{l}E}_{L^2}^2+ \norm{\na^{k+1} B}_{L^2}^2\right)
\\ &\quad\le  C \sum_{l=k}^{k+2}\norm{\na^{l} u}_{L^2}^2
 +C_k  \left(\norm{ (n, u)}_{H^{k+1}\cap H^3}^2+\norm{\na B}_{L^2}^2  \right)  \left(\sum_{l=k}^{k+2}\norm{ \na^{l}(n, u )}_{L^2}^2+\norm{\na^{k+1}B}_{L^2}^2 \right).
\end{split}
\end{equation}
Adding the inequality above with \eqref{n estimate}, we get \eqref{other dissipation}.
\end{proof}

%%%%%%%%%%%%%%%%%%%%%%%%%%%%%%%%%%%%%%%%%%%%%%%%%%%%%%%%%%%%%%%%%%%%%%%%%%%%%%
\subsection{Negative Sobolev estimates}
%%%%%%%%%%%%%%%%%%%%%%%%%%%%%%%%%%%%%%%%%%%%%%%%%%%%%%%%%%%%%%%%%%%%%%%%%%%%%%

In this subsection, we will derive the evolution of the negative Sobolev
norms of $(u,E,B)$. In order to estimate the nonlinear terms, we need to
restrict ourselves to that $s\in (0,3/2)$. We will establish the following lemma.

\begin{lemma}
\label{1Esle} For $s\in(0,1/2]$, we have
\begin{equation}  \label{1E_s}
\frac{d}{dt}\norm{(u,E,B)}_{\dot{H}^{-s}}^2    +\la\norm{u}_{\dot{H}^{-s}}^2
\lesssim \left(\norm{(n,u)}_{H^2}^2+\norm{\na
B}_{H^1}^2 \right)\norm{(u,E,B)}_{\dot{H}^{-s}}+\norm{E}_{H^2}^2;
\end{equation}
and for $s\in(1/2,3/2)$, we have
\begin{equation}  \label{1E_s2}
\begin{split}
&\frac{d}{dt}\norm{(u,E,B)}_{\dot{H}^{-s}}^2  +\la\norm{u}_{\dot{H}^{-s}}^2 \\
&\quad\lesssim\left(\norm{(n,u)}_{H^{1}}^{2}
 +\norm{B}_{L^2}^{s-1/2}\norm{\na B}_{L^2}^{3/2-s}\norm{u}_{L^2}\right)\norm{(u,E,B)}_{\dot{H}^{-s}}+\norm{E}_{H^2}^2.
\end{split}
\end{equation}
\end{lemma}

\begin{proof}
The $\Lambda^{-s}$ $(s>0)$ energy estimate of $\eqref{EM per}_{2}$--$\eqref{EM per}_{4}$ yield
\begin{equation}  \label{1E_s_0}
\begin{split}
&\frac{1}{2}\frac{d}{dt}\norm{(u,E,B)}_{\dot{H}^{-s}}^2    +\nu\norm{u}_{\dot{H}^{-s}}^2
\\&\quad=-\int \Lambda^{-s}\left(u\cdot\na u+\mu  n\na n+u\times B\right)
\cdot\Lambda^{-s} u +\nu \int\Lambda^{-s}(f(n)u)\cdot\Lambda^{-s} E-\int \Lambda^{-s}\na n
\cdot\Lambda^{-s} u
\\
&\quad\lesssim\norm{ u\cdot\na u+\mu  n\na n+u\times B}_{\dot{H}^{-s}} \norm{u}_{\dot{H}^{-s}}
+\norm{f(n)u}_{\dot{H}^{-s}} \norm{E}_{\dot{H}^{-s}} +\norm{ \na n}_{\dot{H}^{-s}} \norm{u}_{\dot{H}^{-s}}.
 \end{split}
\end{equation}

We now restrict the value of $s$ in order to estimate the other terms on
the right-hand side of \eqref{1E_s_0}. If $s\in (0,1/2]$, then $1/2+s/3<1$
and $3/s\geq 6$. Then applying Lemma \ref{Riesz lemma}, together with H\"{o}lder's, Sobolev's and Young's inequalities, we
obtain
\begin{equation}\label{1E_s_1}
\begin{split}
\norm{ u\cdot\na u }
_{\dot{H}^{-s}}
&  \lesssim \norm{ u\cdot\na u }_{L^{\frac{1}{1/2+s/3}}}
 \lesssim \norm{u}_{L^{3/s}}
\norm{\na u }_{L^{2}} \\
&\lesssim \norm{\na u}_{L^{2}}^{1/2+s}
\norm{\na^2u }_{L^{2}}^{1/2-s}\norm{\na u}_{L^2}  \\
& \lesssim \norm{ \na u}_{H^{1}}^{2}+\norm{\na u}_{L^{2}}^{2}  .
\end{split}
\end{equation}
Similarly, we can bound
\begin{eqnarray}
&&\norm{  n\na n }_{\dot{H}^{-s}}\lesssim   \norm{ \na n}
_{H^{1}}^{2}+\norm{\na n}_{L^{2}}^{2} ;  \\
&&\norm{ u\times B}_{\dot{H}^{-s}}
\lesssim \norm{ \na B}
_{H^{1}}^{2}+\norm{u}_{L^{2}}^{2} ;   \\
&&\norm{f(n)u}_{\dot{H}^{-s}}
\lesssim  \norm{ \na u}
_{H^{1}}^{2}+\norm{n}_{L^{2}}^{2} .
\end{eqnarray}

Now if $s\in (1/2,3/2)$, we shall estimate the right-hand side of
\eqref{1E_s_0} in a different way. Since $s\in
(1/2,3/2)$, we have that $1/2+s/3<1$ and $2<3/s<6$. Then applying Lemma \ref{Riesz lemma} and using
(different) Sobolev's inequality, we have
\begin{eqnarray}
&& \nonumber
\norm{ u\cdot\na u }
_{\dot{H}^{-s}} \lesssim \norm{ u}_{L^{3/s}}\norm{\na u }
_{L^{2}} \lesssim \norm{u}_{L^{2}}^{s-1/2}
\norm{\na u }_{L^{2}}^{3/2-s}\norm{\na u}_{L^2}  \\
&&\qquad\qquad  \ \ \ \ \ \lesssim  \norm{u}_{H^{1}}^{2}+\norm{\na u}_{L^{2}}^{2} ;
 \\
&&\norm{  n\na n }_{\dot{H}^{-s}}\lesssim  \norm{n}_{H^{1}}^{2}+\norm{\na n}_{L^{2}}^{2} ;
   \\&&\label{1E_s_52}
\norm{ u\times B}_{\dot{H}^{-s}}
\lesssim \norm{B}
_{L^{2}}^{s-1/2}\norm{ \na  B}_{L^{2}}^{3/2-s}\norm{u}
_{L^{2}};\\
&&\norm{f(n)u}_{\dot{H}^{-s}}
\lesssim  \norm{u}_{H^{1}}^{2}+\norm{n}_{L^{2}}^{2}
 .
\end{eqnarray}

Note that we fail to estimate the remaining last term on the
right-hand side of \eqref{1E_s_0} as above. To overcome this
obstacle, the key point is to make full use of the equation
$\eqref{EM per}_5$ to rewrite $n= n-f(n)+\nu^{-1}{\rm div} E$. This
idea was also used in \cite{TW}. Indeed, using $\eqref{EM per}_5$,
we have
\begin{equation}  \label{na n u}
\begin{split}
\norm{ \na n}_{\dot{H}^{-s}}&\lesssim \norm{\Lambda^{-s}\na {\rm div} E}_{L^2}+\norm{\na(f(n)-n)}_{\dot{H}^{-s}}  \\
&\lesssim  \norm{E}_{H^2}+ \norm{\na(f(n)-n)}_{\dot{H}^{-s}} .
\end{split}
\end{equation}
 Here we have used the facts that $s<3/2$ and $f(n)-n=O(n^2)$. Estimating the last term in \eqref{na n u} as before, and then collecting all the estimates we have derived, by Cauchy's inequality, we deduce \eqref{1E_s} for $s\in(0,1/2]$ and \eqref{1E_s2} for $s\in(1/2,3/2)$.
\end{proof}

%%%%%%%%%%%%%%%%%%%%%%%%%%%%%%%%%%%%%%%%%%%%%%%%%%%%%%%%%%%%%%%%%%%%%%%%%%%%%%
\subsection{Negative Besov estimates}
%%%%%%%%%%%%%%%%%%%%%%%%%%%%%%%%%%%%%%%%%%%%%%%%%%%%%%%%%%%%%%%%%%%%%%%%%%%%%%

In this subsection, we will derive the evolution of the negative Besov
norms of $(u,E,B)$. The argument is similar to the previous subsection.

\begin{lemma}
\label{1Esle2} For $s\in(0,1/2]$, we have
\begin{equation}  \label{1E_s Besov}
\begin{split}
\frac{d}{dt}\norm{(u,E,B)}_{\dot{B}_{2,\infty}^{-s}}^2  +\la\norm{u}_{\dot{B}_{2,\infty}^{-s}}^2 \lesssim \left(\norm{(n,u)}_{H^2}^2+\norm{\na
B}_{H^1}^2 \right)\norm{(u,E,B)}_{\dot{B}_{2,\infty}^{-s}}+\norm{E}_{H^2}^2;
\end{split}
\end{equation}
and for $s\in(1/2,3/2]$, we have
\begin{equation}  \label{1E_s2 Besov}
\begin{split}
&\frac{d}{dt}\norm{(u,E,B)}_{\dot{B}_{2,\infty}^{-s}}^2  +\la\norm{u}_{\dot{B}_{2,\infty}^{-s}}^2 \\
&\quad\lesssim\left(\norm{(n,u)}_{H^{1}}^{2}
 +\norm{B}_{L^2}^{s-1/2}\norm{\na B}_{L^2}^{3/2-s}\norm{u}_{L^2}\right)\norm{(u,E,B)}_{\dot{B}_{2,\infty}^{-s}}+\norm{E}_{H^2}^2.
\end{split}
\end{equation}
\end{lemma}

\begin{proof}
The $\dot{\Delta}_{j}$ energy estimates of $\eqref{EM per}_2$--$\eqref{EM per}_4$ yield, with multiplication of $2^{-2sj}$ and then taking the supremum over ${j\in\mathbb{Z}}$,
\begin{equation}   \label{1E_s_0 Besov}
\begin{split}
\frac{1}{2}&\frac{d}{dt}\norm{(u,E,B)}_{\dot{B}_{2,\infty}^{-s}}^2  +\nu \norm{u}_{\dot{B}_{2,\infty}^{-s}}^2\\
&\lesssim\sup\limits_{j\in\mathbb{Z}}2^{-2sj}\left(  -\int  \dot{\Delta}_{j}\left(u\cdot\na u+\mu  n\na n+u\times B\right)
\cdot\dot{\Delta}_{j} u \right)\\
&\quad+\sup\limits_{j\in\mathbb{Z}}2^{-2sj}\left(\nu \int \dot{\Delta}_{j}(f(n)u)\cdot\dot{\Delta}_{j} E -\int \dot{\Delta}_{j}\na n
\cdot\dot{\Delta}_{j} u\right)
\\
&\lesssim\norm{ u\cdot\na u+\mu  n\na n+u\times B}_{\dot{B}_{2,\infty}^{-s}} \norm{u}_{\dot{B}_{2,\infty}^{-s}}
+\norm{f(n)u}_{\dot{B}_{2,\infty}^{-s}} \norm{E}_{\dot{B}_{2,\infty}^{-s}} +\norm{ \na n}_{\dot{B}_{2,\infty}^{-s}} \norm{u}_{\dot{B}_{2,\infty}^{-s}} .
\end{split}
\end{equation}
Then the proof is exactly the same as the proof of Lemma \ref{1Esle} except that we should apply Lemma \ref{Lp embedding} instead to estimate the $\dot{B}_{2,\infty}^{-s}$ norm. Note that we allow $s=3/2$.
\end{proof}

%%%%%%%%%%%%%%%%%%%%%%%%%%%%%%%%%%%%%%%%%%%%%%%%%%%%%%%%%%%%%%%%%%%%%%%%%%%%%%
\section{Proof of theorems}\label{section3}
%%%%%%%%%%%%%%%%%%%%%%%%%%%%%%%%%%%%%%%%%%%%%%%%%%%%%%%%%%%%%%%%%%%%%%%%%%%%%%

%%%%%%%%%%%%%%%%%%%%%%%%%%%%%%%%%%%%%%%%%%%%%%%%%%%%%%%%%%%%%%%%%%%%%%%%%%%%%%
\subsection{Proof of Theorem \ref{existence}}
%%%%%%%%%%%%%%%%%%%%%%%%%%%%%%%%%%%%%%%%%%%%%%%%%%%%%%%%%%%%%%%%%%%%%%%%%%%%%%

In this subsection, we will prove the unique global solution to the system \eqref{EM per}, and the key point is that we only assume the $H^3$ norm of initial data is small.

{\it Step 1. Global small $\mathcal{E}_3$ solution.}

We first close the energy estimates at the $H^3$ level by assuming a priori that $ \sqrt{\mathcal{E}_3(t)}\le \delta$ is sufficiently small.
Taking $k=0,1$ in \eqref{energy 1} of Lemma \ref{energy lemma} and then summing up, we obtain
\begin{equation} \label{end 1}
\begin{split}
\frac{d}{dt}  \sum_{l=0}^3 \norm{ \na^{l}(n, u, E, B )}_{L^2}^2 +\la\sum_{l=0}^3\norm{\na^{l} u}_{L^2}^2 \lesssim  \sqrt{\mathcal{E}_3}\mathcal{D}_3+\sqrt{\mathcal{D}_3}\sqrt{\mathcal{D}_3}\sqrt{\mathcal{E}_3} \lesssim  \delta\mathcal{D}_3.\end{split}
\end{equation}
Taking $k=0,1$ in \eqref{other dissipation} of Lemma \ref{other di} and then summing up, we obtain
\begin{equation}  \label{end 2}
\begin{split}
&\frac{d}{dt}\left(\sum_{l=0}^2\int \na^lu\cdot\na\na^{l} n +\sum_{l=0}^2\int \na^{l}u\cdot\na^lE  -\eta\sum_{l=0}^1\int
\na^l E \cdot\na^{l}\na\times B \right)
\\&\quad+\la \left(\sum_{l=0}^3\norm{\na^{l}n}_{L^2}^2+\sum_{l=0}^2\norm{\na^{l}E}_{L^2}^2+ \sum_{l=1}^2 \norm{\na^{l} B}_{L^2}^2\right)
\\ &\qquad\lesssim
 \sum_{l=0}^3\norm{\na^{l} u}_{L^2}^2
 +\delta^2\mathcal{D}_3 .
\end{split}
\end{equation}
Multiplying \eqref{end 2} by a sufficiently small but fixed factor $\varepsilon$ and then adding it with \eqref{end 1}, since $\delta$ is small, we deduce that there exists an instant energy functional $\widetilde{\mathcal{E}}_3$ equivalent to ${\mathcal{E}}_3$ such that
\begin{equation}
\frac{d}{dt}\widetilde{\mathcal{E}}_3+\mathcal{D}_3\le 0.
\end{equation}
Integrating the inequality above directly in time, we obtain \eqref{energy inequality}. By a standard continuity argument, we then close the a priori estimates if we assume at initial time that $\mathcal{E}_3(0)\le \delta_0$ is sufficiently small. This concludes the unique global small $\mathcal{E}_3$ solution.

{\it Step 2. Global $\mathcal{E}_N$ solution.}

We shall prove this by an induction on $N\ge 3$. By \eqref{energy inequality}, then \eqref{energy inequality N} is valid for $N=3$. Assume \eqref{energy inequality N} holds for $N-1$ (then now $N\ge 4$). Taking $k=0,\dots,N-2$ in \eqref{energy 1} of Lemma \ref{energy lemma} and then summing up, we obtain
\begin{equation} \label{end 3}
\begin{split}
& \frac{d}{dt}\sum_{l=0}^{N} \norm{ \na^{l}(n,u,E,B)}_{L^2}^2 +\la\sum_{l=0}^{N}\norm{\na^{l} u}_{L^2}^2 \\& \quad\le C_N \sqrt{\mathcal{D}_{N-1}}\sqrt{\mathcal{E}_N}\sqrt{\mathcal{D}_N}+C\sqrt{\mathcal{D}_{N-1}}\sqrt{\mathcal{D}_N}\sqrt{\mathcal{E}_N}
\le C_N \sqrt{\mathcal{D}_{N-1}}\sqrt{\mathcal{E}_N}\sqrt{\mathcal{D}_N}.\end{split}
\end{equation}
Here we have used the fact that $3\le \frac{N-2}{2}+2\le N-2+1$ since $N\ge 4$.
Note that it is important that we have put the two first factors in \eqref{energy 1} into the dissipation.

Taking  $k=0,\dots,N-2$ in \eqref{other dissipation} of Lemma \ref{other di} and then summing up, we obtain
\begin{equation}  \label{end 4}
\begin{split}
&\frac{d}{dt}\left(\sum_{l=0}^{N-1}\int \na^lu\cdot\na\na^{l} n +\sum_{l=0}^{N-1}\int \na^{l}u\cdot\na^lE -\eta\sum_{l=0}^{N-2}\int
\na^l E \cdot\na\times\na^{l} B \right)
\\&\quad+\la \left(\sum_{l=0}^{N}\norm{\na^{l}n}_{L^2}^2+\sum_{l=0}^{N-1}\norm{\na^{l}E}_{L^2}^2+ \sum_{l=1}^{N-1} \norm{\na^{l} B}_{L^2}^2\right)
 \\ &\qquad\le
 C\sum_{l=0}^{N}\norm{\na^{l} u}_{L^2}^2
 +C_N  \sqrt{\mathcal{D}_{N-1}}\sqrt{\mathcal{D}_N}\sqrt{\mathcal{E}_N}.\end{split}
\end{equation}
Multiplying \eqref{end 4} by a sufficiently small factor $\varepsilon$ and then adding it with \eqref{end 3},  we deduce that there exists an instant energy functional $\widetilde{\mathcal{E}}_N$ equivalent to $\mathcal{E}_N$ such that, by Cauchy's inequality,
\begin{equation}
\frac{d}{dt}\widetilde{\mathcal{E}}_N+\mathcal{D}_N\le  C_N \sqrt{\mathcal{D}_{N-1}}\sqrt{\mathcal{E}_N}\sqrt{\mathcal{D}_N}
 \le \tilde{\varepsilon} \mathcal{D}_N+ C_{N,\tilde{\varepsilon}} {\mathcal{D}_{N-1}} {\mathcal{E}_N}
.
\end{equation}
This implies
\begin{equation}\label{EN}
\frac{d}{dt}\widetilde{\mathcal{E}}_N+\frac{1}{2}\mathcal{D}_N
 \le  C_{N } {\mathcal{D}_{N-1}} {\mathcal{E}_N}.
\end{equation}
We then use the standard Gronwall lemma and the induction hypothesis to deduce that
\begin{equation}
\begin{split}
 \mathcal{E}_N(t)+\int_0^t\mathcal{D}_N(\tau)\,d\tau
& \le  C\mathcal{E}_N(0)e^{C_{N }\int_0^t{\mathcal{D}_{N-1}} (\tau)\,d\tau}
\\&\le  C\mathcal{E}_N(0)e^{C_{N }P_{N-1}\left(\mathcal{E}_{N-1}(0)\right)}\\&\le  C\mathcal{E}_N(0)e^{C_{N }P_{N-1}\left(\mathcal{E}_{N}(0)\right)}\equiv P_N\left(\mathcal{E}_{N}(0)\right)
.\end{split}
\end{equation}
This concludes the global $\mathcal{E}_N$ solution. The proof of Theorem \ref{existence} is completed.\hfill$\Box$

\subsection{Proof of Theorem \ref{decay}}
In this subsection, we will prove the various time decay rates of the unique global solution to the system \eqref{EM per} obtained in Theorem \ref{existence}. Fix $N\ge 5$. We need to assume that $\mathcal{E}_N(0)\le \delta_0=\delta_0(N)$ is small. Then Theorem \ref{existence} implies that there exists a unique global $\mathcal{E}_N$ solution, and $\mathcal{E}_N(t)\le P_{N}\left(\mathcal{E}_N(0)\right) \le \delta_0$ is small for all time $t$. Since now our $\delta_0$ is relative small with respect to $N$, we just ignore the $N$ dependence of the constants in the energy estimates in the previous section.

{\it Step 1. Basic decay.}

For the convenience of presentations, we define a family of energy functionals and the corresponding dissipation rates with {\it minimum derivative counts} as
\begin{equation}\label{1111}
\mathcal{E}_{k}^{k+2}=\sum_{l=k}^{k+2}\norm{ \na^{l}(n, u, E, B )}_{L^2}^2
\end{equation}
and
\begin{equation}\label{2222}
\mathcal{D}_{k}^{k+2}=\sum_{l=k}^{k+2}\norm{\na^{l}(n,u)}_{L^2}^2+\sum_{l=k}^{k+1}\norm{\na^{l}E}_{L^2}^2+ \norm{\na^{k+1} B}_{L^2}^2.
\end{equation}

By Lemma \ref{energy lemma}, we have that for $k=0,\dots,N-2$,
\begin{equation}  \label{end 5}
\begin{split}
&\frac{d}{dt}\sum_{l=k}^{k+2}\norm{ \na^{l}(n, u, E, B )}_{L^2}^2 +\la\sum_{l=k}^{k+2}\norm{\na^{l} u}_{L^2}^2 \\
&\quad\lesssim \sqrt{\delta_0} \mathcal{D}_k^{k+2} +\norm{(n,u)}_{L^\infty}\norm{\na^{k+2}(n, u )}_{L^2}
\norm{ \na^{k+2}( E,  B )}_{L^2}.
\end{split}
\end{equation}
By Lemma \ref{other di}, we have that for $k=0,\dots,N-2$,
\begin{equation} \label{end 6}
\begin{split}
&\frac{d}{dt}\left(\sum_{l=k}^{k+1}\int \na^lu\cdot\na\na^{l} n +\sum_{l=k}^{k+1}\int \na^{l}u\cdot\na^lE  -\eta\int
\na^k E \cdot\na^{k}\na\times B \right)
\\&\quad+\la \left(\sum_{l=k}^{k+2}\norm{\na^{l}n}_{L^2}^2+\sum_{l=k}^{k+1}\norm{\na^{l}E}_{L^2}^2+ \norm{\na^{k+1} B}_{L^2}^2\right)
\\ &\qquad\lesssim\sum_{l=k}^{k+2}\norm{\na^{l} u}_{L^2}^2
 +\delta_0 \sum_{l=k}^{k+2}\norm{\na^{l} (n, u )}_{L^2}^2.
\end{split}
\end{equation}
Multiplying \eqref{end 6} by a sufficiently small but fixed factor $\varepsilon$ and then adding it with \eqref{end 5}, since $\delta_0$ is small, we deduce that there exists an instant energy functional $\widetilde{\mathcal{E}}_k^{k+2}$ equivalent to $\mathcal{E}_k^{k+2}$ such that
\begin{equation}\label{energy}
\frac{d}{dt}\widetilde{\mathcal{E}}_k^{k+2}+\mathcal{D}_k^{k+2}\lesssim   \norm{(n,u)}_{L^\infty}\norm{\na^{k+2}(n, u )}_{L^2}
\norm{ \na^{k+2}( E,  B )}_{L^2}.
\end{equation}
Note that we cannot absorb the right-hand side of \eqref{energy} by the dissipation $\mathcal{D}_k^{k+2}$ since it does not contain $\norm{\na^{k+2}( E, B )}_{L^2}^2$. We will distinguish the arguments by the value of $k$. If $k=0$ or $k=1$, we bound $\norm{\na^{k+2}( E, B )}_{L^2}$ by the energy. Then we have that for $k=0,1$,
\begin{equation}
\frac{d}{dt}\widetilde{\mathcal{E}}_k^{k+2}+\mathcal{D}_k^{k+2}\lesssim   \sqrt{\mathcal{D}_k^{k+2}}\sqrt{\mathcal{D}_k^{k+2}}\sqrt{{\mathcal{E}}_3}\lesssim \sqrt{\delta_0}\mathcal{D}_k^{k+2},
\end{equation}
which implies
\begin{equation}
\frac{d}{dt}\widetilde{\mathcal{E}}_k^{k+2}+\mathcal{D}_k^{k+2}\le 0.
\end{equation}
If $k\ge 2$, we have to bound $\norm{\na^{k+2}( E, B )}_{L^2}$ in term of $\norm{\na^{k+1}( E, B )}_{L^2}$ since $\sqrt{\mathcal{D}_k^{k+2}}$ cannot control $\norm{(n,u)}_{L^\infty}$. The key point is to use the regularity interpolation method developed in \cite{GW,SG06}. By Lemma \ref{A1}, we have
\begin{equation}\label{kkk}
\begin{split}
&\norm{(n,u)}_{L^\infty}\norm{ \na^{k+2}(n,u )}_{L^2} \norm{ \na^{k+2}( E, B )}_{L^2}\\&\quad\lesssim\norm{(n,u)}_{L^2}^{1-\frac{3 }{2k}}\norm{\na^{k}(n,u )}_{L^2}^{\frac{3 }{2k}}\norm{\na^{k+2}(n,u )}_{L^2}\norm{\na^{k+1}( E, B )}_{L^2}^{1-\frac{3 }{2k}}
\norm{\na^{\al}( E, B )}_{L^2}^{\frac{3 }{2k}},
\end{split}
\end{equation}
where $\alpha$ is defined by
\begin{equation}
 k+2=(k+1)\times\left(1-\frac{3 }{2k}\right)+\alpha\times \frac{3 }{2k}
 \Longrightarrow \alpha=\frac{5}{3}k+1.
\end{equation}
Hence, for $k\ge 2$, if $N\ge \frac{5}{3}k+1\Longleftrightarrow 2\le k\le \frac{3}{5}(N-1)$, then by \eqref{kkk}, we deduce from \eqref{energy} that
\begin{equation}
\frac{d}{dt}\widetilde{\mathcal{E}}_k^{k+2}+\mathcal{D}_k^{k+2} \lesssim  \sqrt{{\mathcal{E}}_{N}} {\mathcal{D}_k^{k+2}}\lesssim \sqrt{\delta_0}\mathcal{D}_k^{k+2},
\end{equation}
which allow us to arrive at that for any integer $k$ with $0\le k\le  \frac{3}{5}(N-1) $ (note that $N-2\ge \frac{3}{5}(N-1)\ge 2$ since $N\ge 5$), we have
\begin{equation}\label{k energy}
\frac{d}{dt}\widetilde{\mathcal{E}}_k^{k+2}+\mathcal{D}_k^{k+2} \le 0.
\end{equation}

The fact that $\mathcal{D}_k^{k+2}$ is weaker than $\mathcal{E}_k^{k+2}$ prevents the exponential decay of the solution. In order to effectively derive the decay rate from \eqref{k energy}, we still manage to bound the missing terms in the energy, that is, $\norm{\na^k B}_{L^2}^2$ and $\norm{ \na^{k+2}( E, B )}_{L^2}^2$ in terms of $\mathcal{E}_k^{k+2}$ in \eqref{k energy}. We again use the regularity interpolation method, but now we need to also do the Sobolev interpolation between the negative and positive Sobolev norms.  Assuming for the moment that we
have proved \eqref{H-sbound} or \eqref{H-sbound Besov}. Using Lemma
\ref{1-sinte}, we have that for $s\ge 0$ and $k+s>0$,
\begin{equation}\label{inter 1'}
\norm{\na^k B }_{L^2} \le  \norm{B}_{\dot{H}^{-s}}^{\frac{1}{k+1+s}}\norm{\na^{k+1}B }_{L^2}^{\frac{k+s}{k+1+s}}
 \le  C_0\norm{\na^{k+1}B}_{L^2}^{\frac{k+s}{k+1+s}}.
\end{equation}
Similarly, using Lemma
\ref{Besov interpolation}, we  have that for $s>0$ and $k+s>0$,
\begin{equation}\label{inter 11'}
\norm{\na^k B }_{L^2} \le  \norm{B }_{\dot{B}_{2,\infty}^{-s}}^{\frac{1}{k+1+s}}\norm{\na^{k+1}B }_{L^2}^{\frac{k+s}{k+1+s}}
 \le  C_0\norm{\na^{k+1}B}_{L^2}^{\frac{k+s}{k+1+s}}.
\end{equation}
On the other hand, for $k+2<N$, we have
\begin{equation}
\begin{split}
\norm{ \na^{k+2}( E, B )}_{L^2}\le  \norm{\na^{k+1}( E, B )}_{L^2}^{ \frac{N-k-2}{N-k-1}}
\norm{\na^N( E, B )}_{L^2}^{\frac{1}{N-k-1}}\le  C_0\norm{\na^{k+1}( E, B )}_{L^2}^{ \frac{N-k-2}{N-k-1}}.
\end{split}
\end{equation}
Then we deduce from \eqref{k energy} that
\begin{equation}
\frac{d}{dt}\widetilde{\mathcal{E}}_k^{k+2}+\left\{\mathcal{E}_k^{k+2}\right\}^{1+\vartheta} \le 0,
\end{equation}
where $\vartheta=\max\left\{\frac{1}{k+s},\frac{1}{N-k-2}\right\}$.
Solving this inequality directly, we obtain in particular that
\begin{equation}\label{nnn}
\mathcal{E}_k^{k+2}  (t) \le \left\{\left[\mathcal{E}_k^{k+2}(0)\right]^{-\vartheta}+\vartheta  t\right\}^{-  {1}/{\vartheta}}
\le C_0 (1+ t)^{- {1}/{\vartheta}}=C_0 (1+ t)^{-\min\left\{ {k+s}, {N-k-2}\right\}}.
\end{equation}
Notice that \eqref{nnn} holds also for $k+s=0$ or $k+2=N$. So, if we want to obtain the optimal decay rate of the whole solution for the spatial derivatives of order $k$, we only need to assume $N$ large enough (for fixed $k$ and $s$) so  that $k+s\le N-k-2$. Thus we should require that
\begin{equation}
N\ge \max\left\{k+2, \frac{5}{3}k+1, 2k+2+s\right\}= 2k+2+s.
\end{equation}
This proves the optimal decay \eqref{basic decay}.

Finally, we turn back to prove \eqref{H-sbound} and \eqref{H-sbound Besov}. First, we prove \eqref{H-sbound} by using Lemma \ref{1Esle}. However, we are not
able to prove them for all $s\in[0,3/2)$ at this moment. We must distinguish the arguments by the value of $s$.
First, for $s\in (0,1/2]$, integrating \eqref{1E_s} in time, by \eqref{energy inequality} we obtain that for $s\in (0,1/2]$,
\begin{equation}
\begin{split}
\norm{(u,E,B)(t)}_{\dot{H}^{-s}}^2& \lesssim \norm{(u_0,E_0,B_0) }_{\dot{H}^{-s}}^2 + \int_{0}^{t} \mathcal{D}_3(\tau ) \left(1+\norm{(u,E,B)(\tau)}_{\dot{H}^{-s}}\right)
 \,d\tau  \\
& \leq  C_0\left(1+\sup_{0\leq \tau \leq t}\norm{(u,E,B)(\tau)}_{\dot{H}^{-s}}\right).
\end{split}
\label{1-sin2}
\end{equation}
By Cauchy's inequality, this together with \eqref{energy inequality}
gives \eqref{H-sbound} for $s\in [ 0,1/2]$ and thus verifies \eqref{basic decay} for $s\in [ 0,1/2]$. Next, we let $s\in (1/2,1)$.
Observing that we have $ (u _0,E_0,B_0)\in \dot{H}^{-1/2}
$ since $\dot{H}^{-s}\cap L^2\subset\dot{H}^{-s^{\prime}}$ for any $
s^{\prime}\in [0,s]$, we then deduce from what we have proved for
\eqref{basic decay} with $s=1/2$ that the following decay
result holds:
\begin{equation}  \label{1proof14}
\norm{\na^k  (n,u,E,B)(t)}
_{L^{2}}  \le C_0(1+t)^{-\frac{k+  {1}/{2}}{2}}\ \hbox{ for }k=0,1.
\end{equation}
Here, since we have required $N\ge 5$ and now $s=1/2$, we have used $k=1$ in \eqref{basic decay}. Thus by \eqref{1proof14}, \eqref{energy inequality} and H\"older's inequality, we deduce from \eqref{1E_s2} that for $
s\in(1/2,1)$,
\begin{equation}  \label{1-sin2''}
\begin{split}
\norm{(u,E,B)(t)}_{\dot{H}^{-s}}^2
& \lesssim \norm{(u_0,E_0,B_0)}_{\dot{H}^{-s}}^2 + \int_{0}^{t} \mathcal{D}_3(\tau ) \left(1+\norm{(u,E,B)(\tau)}_{\dot{H}^{-s}}\right)
 \,d\tau\\
&\quad+\int_0^t\norm{B(\tau )}_{L^2}^{s-1/2}\norm{\na B(\tau )}_{L^2}^{3/2-s}\sqrt{\mathcal{D}_3(\tau )}\norm{(u,E,B)(\tau)}_{\dot{H}^{-s}}\,d\tau\\
&\le C_0\left(1+\left(1+\int_0^t(1+\tau)^{-2(1-s/2)}\,d\tau \right)\sup_{0\le\tau\le t}\norm{(u,E,B)(\tau)}_{\dot{H}^{-s}}\right)\\
 &\le C_0\left(1+\sup_{0\le\tau\le t}\norm{(u,E,B)(\tau)}_{\dot{H}^{-s}}\right).
\end{split}
\end{equation}
Here we have used the fact $s\in(1/2,1)$ so that the time integral in \eqref{1-sin2''} is finite.
This gives \eqref{H-sbound} for $s\in (1/2,1)$ and thus verifies \eqref{basic decay} for $s\in (1/2,1)$.
Now let $s\in [1,3/2)$. We choose $s_0$ so that $s-1/2<s_0<1$. Hence, $ (u _0,E_0,B_0)\in \dot{H}^{-s_0}$. We then deduce from what we have proved for \eqref{basic decay} with $s=s_0$ that the following decay
result holds:
\begin{equation}  \label{2proof14}
\norm{\na^k  (n,u,E,B)(t)}
_{L^{2}}  \le C_0(1+t)^{- \frac{k+s_0}{2} }\ \hbox{ for }k=0,1.
\end{equation}
Here, since we have required $N\ge 5$ and now $s=s_0<1$, we have used $k=1$ in \eqref{basic decay}. Thus by \eqref{2proof14} and H\"older's inequality, we deduce from \eqref{1E_s2} that for $
s\in[1,3/2)$, similarly as in \eqref{1-sin2''},
\begin{equation}\label{mmm}
\begin{split}
 \norm{(u,E,B)(t)}_{\dot{H}^{-s}}^2
 &\le C_0\left(1+\left(1+\int_0^t(1+\tau)^{- ( {s_0} + 3/2-  s )}\,d\tau \right)\sup_{0\le\tau\le t}\norm{(u,E,B)(\tau)}_{\dot{H}^{-s}}\right)\\
 &\le C_0\left(1+\sup_{0\le\tau\le t}\norm{(u,E,B)(\tau)}_{\dot{H}^{-s}}\right).
\end{split}
\end{equation}
Here we have used the fact $s-s_0<1/2$ so that the time integral in \eqref{mmm} is finite.
This gives \eqref{H-sbound} for $s\in [1,3/2)$ and thus verifies \eqref{basic decay} for $s\in [1,3/2)$. Note that
\eqref{H-sbound Besov} can be proved similarly except that we use instead  Lemma \ref{1Esle2}.

{\it Step 2. Further decay.}

We first prove \eqref{further decay1} and \eqref{further decay11}. First, noticing that $-\nu f(n)={\rm div}E$, by \eqref{basic decay} and Lemma \ref{A2}, if $N\ge2k+4+s$, then
\begin{equation}\label{n further decay}
\norm{\na^kn(t)}_{L^2}\lesssim\norm{\na^kf(n)(t)}_{L^2}\lesssim\norm{\na^{k+1}E(t)}_{L^2}\le C_0(1+t)^{-\frac{k+1+s}{2}},
\end{equation}
where we have used $n=f^{-1}(f(n))$.

Next, applying $\na^k$ to $\eqref{EM per}_2, \eqref{EM per}_3$ and then multiplying the resulting identities
by $\na^ku$, $\na^kE$ respectively, summing up and integrating over $\mathbb{R}^3$, we obtain
\begin{equation}\label{uE yi}
\begin{split}
& \frac{1}{2}\frac{d}{dt}\int\norms{\na^k(u,E)}^2+\nu \norm{\na^ku}_{L^2}^2\\
&\quad=-\int\na^k\left(\na n+u\cdot\na u
+\mu n\na n+u\times B\right)\cdot\na^ku+\nu\int\na^k\left( \na\times B+ f(n)u\right)\cdot\na^kE
\\
 &\quad\lesssim \norm{\na^{k+1} n}_{L^2} \norm{\na^{k}  u}_{L^2}  +\norm{\na^{k}\left(u\cdot\na
u+\mu n\na n+u\times B\right)}_{L^2}\norm{\na^{k} u}_{L^2}
\\
 &\qquad+\norm{ \na^{k}\left(\na\times B+f(n)u
\right)}_{L^2}\norm{\na^{k} E}_{L^2}.
\end{split}
\end{equation}
On the other hand, taking $l=k$ in \eqref{identy 1}, we may have
\begin{equation} \label{uE san}
\begin{split}
&\int  \na^k \partial_tu \cdot\na^{k}E +\nu \norm{\na^{k} E}_{L^2}^2
\\
 &\quad\lesssim\left(\norm{\na^{k+1}n}_{L^2}+\norm{\na^{k} u}_{L^2}\right)\norm{\na^{k}  E}_{L^2} +\norm{\na^{k}\left(u\cdot\na
u+\mu n\na n+u\times B\right)}_{L^2}\norm{\na^{k} E}_{L^2}.
\end{split}
\end{equation}
Substituting \eqref{inter 21} with $l=k$ into \eqref{uE san}, we may then have
\begin{equation}\label{uE wu}
\begin{split}
\frac{d}{dt}&\int  \na^ku \cdot\na^{k}E +\nu \norm{\na^{k} E}_{L^2}^2
\\
 &\le C \norm{\na^{k} u}_{L^2}^2+C\left(\norm{\na^{k+1} n}_{L^2}+ \norm{\na^{k} u}_{L^2}\right)\norm{\na^{k}  E}_{L^2} \\
 &\quad+\norm{\na^{k}\left(u\cdot\na
u+\mu n\na n+u\times B\right)}_{L^2}\norm{\na^{k} E}_{L^2}+\norm{ \na^{k}\left(\na\times B+f(n)u
\right)}_{L^2}\norm{\na^{k} u}_{L^2}.
\end{split}
\end{equation}
Multiplying \eqref{uE wu} by a sufficiently small but fixed factor $\varepsilon$ and then adding it with \eqref{uE yi}, since $\varepsilon$ is small, we deduce that there exists $\mathcal{F}_k(t)$ equivalent to $\norm{\na^k(u,E)(t)}_{L^2}^2$ such that, by Cauchy's inequality, \eqref{I1 I2}, \eqref{I3 k}, \eqref{I4 k k+1}, \eqref{basic decay} and \eqref{n further decay},
\begin{align}\label{uE liu}
 &\frac{d}{dt}\mathcal{F}_k(t)+\mathcal{F}_k(t)\nonumber
 \\&\quad\lesssim \norm{\na^{k} n}_{L^2}^2+\norm{\na^{k+1}(n,B)}_{L^2}^2+\norm{\na^{k}\left(u\cdot\na
u+\mu n\na n+u\times B\right)}_{L^2}^2 +\norm{ \na^{k}(f(n)u)
}_{L^2}^2\nonumber
  \\&\quad\lesssim \norm{\na^{k} n}_{L^2}^2+\norm{\na^{k+1} (n, B)}_{L^2}^2 + \left(\norm{u}_{H^{\frac k2}}+\norm{\na B}_{L^2}\right)^2\norm{\na^{k+1}B}_{L^2}^2\nonumber\\
&\qquad+\norm{(n,u)}_{L^\infty}^2\norm{\na^{k+1}(n,u)}_{L^2}^2\nonumber
\\
&\quad\le C_0(1+t)^{-(k+1+s)},
\end{align}
where we required $N\ge2k+4+s$.
Applying the standard Gronwall lemma to \eqref{uE liu}, we obtain
\begin{equation}\label{}
\begin{split}
\mathcal{F}_k(t)\le \mathcal{F}_k(0)e^{-t}+C_0\int_0^te^{-(t-\tau)}(1+\tau)^{-(k+1+s)}\,d\tau\le C_0(1+t)^{-(k+1+s)}.
\end{split}
\end{equation}
This implies
\begin{equation}\label{uE qi}
\norm{\na^k(u,E)(t)}_{L^2}\lesssim \sqrt{\mathcal{F}_k(t)}\le C_0(1+t)^{-\frac{k+1+s}{2}}.
\end{equation}
We thus complete the proof of \eqref{further decay1}. Notice that \eqref{further decay11} now follows by \eqref{n further decay} with the improved decay rate of $E$ in \eqref{further decay1}, just requiring $N\ge2k+6+s$.

Now we prove \eqref{further decay2}. Assuming $ B_\infty =0$, then we can extract the following system from  $\eqref{EM per}_1$--$\eqref{EM per}_2$, denoting $\psi={\rm div} u$,
\begin{equation}\label{npsi yi}
\left\{
\begin{array}{lll}
\displaystyle\partial_tn+\psi=-u\cdot\na n
-\mu n{\rm div} u,   \\
\displaystyle\partial_t \psi+\nu  \psi-\nu^2 n=-\Delta n-{\rm div}(u\cdot\na u
+\mu n\na n+u\times B)+\nu^2(f(n)-n).
\end{array}
\right.
\end{equation}
Applying $\na^k$ to $\eqref{npsi yi}$ and then multiplying the resulting identities
by $\nu^2\na^kn$, $\na^k\psi$, respectively, summing up and integrating over $\mathbb{R}^3$, we obtain
\begin{equation}\label{npsi er}
\begin{split}
\frac{1}{2}&\frac{d}{dt}\int \nu^2\norms{\na^k n }^2+\norms{\na^k \psi }^2+\nu\norm{\na^k\psi}_{L^2}^2\\
&=-\nu^2\int\na^k(u\cdot\na n+\mu n{\rm div} u) \na^kn-\int\na^k \Delta n \na^k\psi\\
&\quad-\int\na^k\left[ {\rm div}(u\cdot\na u
+\mu n\na n+u\times B)-\nu^2(f(n)-n)\right] \na^k\psi.
\end{split}
\end{equation}
Applying $\na^k$ to $\eqref{npsi yi}_2$ and then multiplying by $-\na^kn$, as before integrating by parts over $t$ and $x$ variables and using the equation $\eqref{npsi yi}_1$, we may obtain
%-\int&\na^k \pa_t \psi\cdot\na^kn-\nu\int\na^k \psi\cdot\na^kn+\nu^2\norm{\na^kn}_{L^2}^2\\
%&=\int\na^k[\Delta n+{\rm div}(u\cdot\na u
%+\mu n\na n+u\times B)-\nu^2(f(n)-n)]\cdot\na^kn
\begin{equation}\label{npsi san}
\begin{split}
-&\frac{d}{dt}\int\na^k\psi \na^kn+\nu^2\norm{\na^kn}_{L^2}^2\\
&=\norm{\na^k\psi}_{L^2}^2+\nu\int\na^kn \na^k\psi+\int\na^k(u\cdot\na n+\mu n{\rm div} u) \na^k\psi\\
&\quad+\int\na^k\left[\Delta n+{\rm div}(u\cdot\na u
+\mu n\na n+u\times B)-\nu^2(f(n)-n)\right] \na^kn.
\end{split}
\end{equation}

Multiplying \eqref{npsi san} by a sufficiently small but fixed factor $\varepsilon$ and then adding it with \eqref{npsi er}, since $\varepsilon$ is small, we deduce that there exists $\mathcal{G}_k(t)$ equivalent to $\norm{\na^k(n,\psi)}_{L^2}^2$ such that, by Cauchy's inequality,
\begin{equation}\label{npsi si}
\begin{split}
 \frac{d}{dt}\mathcal{G}_k(t)+\mathcal{G}_k(t)&\lesssim  \norm{\na^{k+2}n}_{L^2}^2+\norm{\na^{k+1}(u\cdot\na u)}_{L^2}^2+\norm{\na^{k+1}(n\na n)}_{L^2}^2 \\
 &\quad +\norm{\na^{k+1}(u\times B)}_{L^2}^2+ \norm{\na^k(u\cdot\na n)}_{L^2}^2+\norm{\na^k(n{\rm div}u)}_{L^2}^2.
 \end{split}
\end{equation}
Notice that we have used the following estimates, by using the arguments in Lemma \ref{A2},
\begin{equation}
\norm{\na^{k}(f(n)-n)}_{L^2}^2\le C_k\norm{n}_{H^3}^2\norm{\na^kn}_{L^2}^2\lesssim\delta_0\norm{\na^kn}_{L^2}^2.
\end{equation}
By Lemma \ref{commutator} and Cauchy's inequality, we obtain
\begin{equation}\label{npsi liu}
\begin{split}
\norm{\na^{k+1}(u\times B)}_{L^2}^2&=\norm{u\times \na^{k+1}B+\left[\na^{k+1}, u\right]\times B}_{L^2}^2\\
&\lesssim\norm{u\times \na^{k+1}B}_{L^2}^2+\norm{\left[\na^{k+1}, u\right]\times B}_{L^2}^2\\
&\lesssim\norm{u}_{L^\infty}^2\norm{\na^{k+1}B}_{L^2}^2+\norm{\na u}_{L^\infty}^2\norm{\na^{k}B}_{L^2}^2+\norm{\na^{k+1}u}_{L^2}^2\norm{B}_{L^\infty}^2.
\end{split}
\end{equation}
The other nonlinear terms on the right-hand side of \eqref{npsi si} can be estimated similarly.
Hence, we deduce from \eqref{npsi si} that, by \eqref{basic decay}--\eqref{further decay11},
\begin{equation}\label{npsi qi}
\begin{split}
 &\frac{d}{dt}\mathcal{G}_k(t)+\mathcal{G}_k(t)
 \\&\quad\lesssim \norm{\na^{k+2}n}_{L^2}^2+\norm{u}_{L^\infty}^2\norm{\na^{k+1}B}_{L^2}^2+\norm{\na u}_{L^\infty}^2\norm{\na^{k}B}_{L^2}^2+\norm{B}_{L^\infty}^2\norm{\na^{k+1}u}_{L^2}^2\\
&\qquad+\norm{(n,u)}_{L^\infty}^2\norm{\na^{k+2}(n,u)}_{L^2}^2 +\norm{\na(n,u)}_{L^\infty}^2\norm{\na^{k+1}(n,u)}_{L^2}^2
\\&\quad\le C_0\left((1+t)^{-(k+3+s)}+  (1+t)^{-(k+7/2+2s)} + (1+t)^{-(k+11/2+2s)} \right)
\\&\quad\le C_0(1+t)^{-(k+3+s)},
 \end{split}
\end{equation}
where we required $N\ge2k+8+s$.
Applying the Gronwall lemma to \eqref{npsi qi} again, we obtain
\begin{equation}
\mathcal{G}_k(t)\le \mathcal{G}_k(0)e^{-t}+C_0\int_0^te^{-(t-\tau)}(1+\tau)^{-(k+3+s)}\,d\tau\le  C_0(1+t)^{-(k+3+s)}.
\end{equation}
This implies
\begin{equation}\label{npsi ba}
\norm{\na^k(n,\psi)(t)}_{L^2}\lesssim \sqrt{\mathcal{G}_k(t) }\le C_0 (1+t)^{-\frac{k+3+s}{2}}.
\end{equation}
If required $N\ge2k+12+s$, then by \eqref{npsi ba}, we have
\begin{equation}\label{npsi jiu}
\norm{\na^{k+2}n(t)}_{L^2} \le C_0(1+t)^{-\frac{k+5+s}{2}}.
\end{equation}
Having obtained such faster decay, we can then improve  \eqref{npsi qi} to be
\begin{equation}\label{npsi shi}
 \frac{d}{dt}\mathcal{G}_k(t)+\mathcal{G}_k(t)
 \le C_0\left((1+t)^{-(k+5+s)}+  (1+t)^{-(k+7/2+2s)}  \right)
  \le C_0(1+t)^{-(k+7/2+2s)}.
\end{equation}
Applying the Gronwall lemma again, we obtain
\begin{equation}\label{npsi shiyi}
\norm{\na^k(n,\psi)(t)}_{L^2}\lesssim \sqrt{\mathcal{G}_k(t) }\le C_0 (1+t)^{-(k/2+7/4+s)}.
\end{equation}
We thus complete the proof of \eqref{further decay2}. The proof of Theorem \ref{decay} is completed.

\end{document}